\documentclass[11pt]{amsart}

\usepackage[top=1in, left=1in, right=1in, bottom=1in]{geometry}

\usepackage{amsmath,amssymb,amsthm,amsfonts}
\usepackage{paralist}
\usepackage{booktabs}
\usepackage[usenames,dvipsnames]{color}
\usepackage{comment}
\usepackage{graphicx,epsf}
\usepackage{graphics}
\usepackage{epsfig}
\usepackage{epstopdf}
\usepackage{psfrag}

\usepackage{amsmath,amsfonts,amssymb,latexsym,amsthm}
\usepackage{graphicx,epsf}
\usepackage{amsthm}
\usepackage{multicol}
\usepackage{ytableau}

\usepackage{url}
\usepackage{hyperref} 

\usepackage{tikz}
\usetikzlibrary{automata,positioning,arrows}

\usepackage[colorinlistoftodos]{todonotes}

\RequirePackage{cleveref}
\usepackage{hypcap}
\hypersetup{colorlinks=true, citecolor=darkblue, linkcolor=darkblue}
\definecolor{darkblue}{rgb}{0.0,0,0.7}
\newcommand{\darkblue}{\color{darkblue}}

\definecolor{darkred}{rgb}{0.68,0,0}

\definecolor{darkgreen}{rgb}{0,.38,0}

\newcommand{\defn}[1]{\emph{\darkblue #1}}

\newlength{\ml}
\newtheorem{theorem}{Theorem}
\newtheorem{thm}{Theorem}[section]

\newtheorem{lemma}[thm]{Lemma}

\newtheorem{problem}[thm]{Problem}

\newtheorem{cor}[thm]{Corollary}
\newtheorem{prop}[thm]{Proposition}

\newtheorem{conjecture}[thm]{Conjecture}
\newtheorem{conj}[thm]{Conjecture}

\newtheorem{op}[thm]{Open Problem}
\newtheorem{ex}[thm]{Example}

\theoremstyle{plain}

\theoremstyle{definition}

\newtheorem{rem}[thm]{Remark}

\numberwithin{equation}{section}

\def\bx{{\textbf{\textit{x}}}}
\def\by{{\textbf{\textit{y}}}}

\newcommand{\la}{\lambda}

\newcommand{\cc}{\Bbb C}

\newcommand{\ov}{\overline}

\newcommand{\ga}{\gamma}
\newcommand{\si}{\sigma}
\newcommand{\de}{\delta}

\newcommand{\al}{\alpha}
\newcommand{\be}{\beta}

\newcommand{\A}{{\mathbb A}}

\def\per{\mathrm{per}}
\def\sgn{\mathrm{sgn}}
\def\sk{{\textrm{sk}}}

\def\.{\hskip.06cm}
\def\ts{\hskip.03cm}

\def\SSYT{ {\text {\rm SSYT}  } }

\def\GL{ {\text {\rm GL}  } }

\def\P{{{\rm{\textsf{P}} }}}

\def\FP{{\rm{\textsf{FP}}}}
\def\PH{{\rm{\textsf{PH}}}}

\def\SP{{{\rm{\textsf{\#P}}}}}

\newcommand{\ComCla}[1]{\textup{\textsf{#1}} }
\def\BQP{\ComCla{BQP}}
\def\QMA{\ComCla{QMA}}

\def\GapP{{{\rm{\textsf{GapP}} }}}
\newcommand{\CeqP}{\ComCla{C$_=$P}}
\def\FP{{{}\rm{\textsf{FP}} }}
\def\NP{{{\rm{\textsf{NP}} }}}
\def\coNP{{{\rm{\textsf{coNP}} }}}

\def\VP{{{\rm{\textsf{VP}} }}}
\def\VBP{{{\rm{\textsf{VBP}} }}}
\def\VNP{{{\rm{\textsf{VNP}} }}}

\def\N{\mathbb{N}}

\def\C{\mathbb{C}}
\def\IC{{\mathbb{C}}}

\def\VPs{\mathsf{VP}_{\mathrm{ws}}}
\def\dc{\mathsf{dc}}

\def\Pow{\mathsf{Pow}}

\def\<{\langle}
\def\>{\rangle}

\def\N{\mathbb{N}}

\def\C{\mathbb{C}}

\def\det{\mathrm{det}}
\def\per{\mathrm{per}}
\def\sgn{\mathrm{sgn}}
\def\GL{\mathrm{GL}}

\def\Sym{\mathsf{Sym}}
\newcommand{\Det}{\Omega} 
\newcommand{\ol}[1]{\overline{#1}}

\def\la{\lambda}

\def\mult{\mathrm{mult}}

\def\id{\mathrm{id}}

\def\ocp{Z}

\newcommand{\IA}{\ensuremath{\mathbb{A}}}

\newcommand{\Ch}{\textup{\textsf{Ch}}}

\title[CC in AC]{Computational Complexity in Algebraic Combinatorics \\[2ex] \small{ Current Developments in Mathematics 2023\\Talk notes} }

\author{Greta Panova}

\address{Department of Mathematics, University of Southern California, 
Los Angeles, CA 90089}
\email{gpanova@usc.edu}
\urladdr{https://sites.google.com/usc.edu/gpanova/}

\begin{document}
\maketitle

{\hfill {\footnotesize ``It's hard to look for a black cat in a dark room, especially if there is no cat...''}}
 
{\hfill {\footnotesize Confucius}}

\begin{abstract}

Algebraic Combinatorics originated in Algebra and Representation Theory, studying their discrete objects and integral quantities via combinatorial methods which have since developed independent and self-contained lives and brought us some beautiful formulas and combinatorial interpretations. 
The flagship hook-length formula counts the number of Standard Young Tableaux, which also gives the dimension of the irreducible Specht modules of the Symmetric group. The elegant Littlewood-Richardson rule gives the multiplicities of irreducible GL-modules in the tensor products of GL-modules. Such formulas and rules have inspired large areas of study and development beyond Algebra and Combinatorics, becoming applicable to Integrable Probability and Statistical Mechanics, and Computational Complexity Theory. 

We will see what lies beyond the reach of such nice product formulas and combinatorial interpretations and enter the realm of Computational Complexity Theory, that could formally explain the beauty we see and the difficulties we encounter in finding further formulas and ``combinatorial interpretations''. A 85-year-old such problem asks for a positive combinatorial formula for the Kronecker coefficients of the Symmetric group, another one pertains to the plethysm coefficients of the General Linear group. 

In the opposite direction,  the study of Kronecker and plethysm coefficients leads to the disproof of the wishful approach of Geometric Complexity Theory (GCT) towards the resolution of the algebraic P vs NP Millennium problem, the VP vs VNP problem. In order to make GCT work and establish computational complexity lower bounds, we need to understand representation theoretic multiplicities in further detail, possibly asymptotically. 

\end{abstract}

\section{Introduction}

It was the best of times, it was the worst of times, it was the epoch of unwavering faith in mathematical conjectures, it was the era of crushing despair in the wake of disproofs. In the realm of Algebraic Combinatorics, where beauty had long flourished in the form of graceful formulas and elegant combinatorial interpretations, a shifting tide of uncanny difficulties now swept across the landscape. The hope for solely aesthetic solutions would fade in the shadows of the rigorous framework of Computational Complexity Theory and the imprecision of asymptotic analysis.

What is  Algebraic Combinatorics? According to Wikipedia, it ``is an area of mathematics that employs methods of abstract algebra, notably group theory and representation theory, in various combinatorial contexts and, conversely, applies combinatorial techniques to problems in algebra.'' Here, we will narrow it down to the intersection of Representation Theory and Discrete Mathematics, and in more concrete terms the area consisting of Symmetric Function Theory, and Representation Theory of $S_n$ and $GL_N$ which house our favorite standard and semi-standard Young tableaux. 

The counterpart in our study, Computational Complexity theory, is about the classification of computational problems by how efficiently with respect to the given resource (space or time) they can be solved by an algorithm. It is home to the $\P$ vs $\NP$ problem, and its algebraic version the $\VP$ vs $\VNP$ problem. 

The two fields come together in two ways. Algebraic Combinatorics has old classical quantities (structure constants) resisting a formula or even a ``combinatorial interpretation'' for more than 80 years. Computational Complexity can formalize these difficulties and explain what is [not] happening. On the other side, these same structure constants appear central to the problem of showing that $\VP \neq \VNP$ using Geometric Complexity Theory, and in particular the search for ``multiplicity obstructions''. 

\subsection{Problems in Algebraic Combinatorics through the prism of Computational Complexity}

The dawn of Algebraic Combinatorics was lit by beautiful formulas and elegant combinatorial interpretations. The number $f^\la$ of standard Young tableaux of shape $\la$ is given by the hook-length formula of Frame, Robinson and Thrall. From the representation theoretic correspondences, this is also the dimension of the irreducible $S_n$ representation, the Specht module $\mathbb{S}_\la$.  The tensor product of two $GL_N$ irreducibles $V_\mu$ and $V_\nu$ factors into irreducibles with multiplicities given by the Littlewood-Richardson coefficients $c^\la_{\mu\nu}$.  While no ``nice'' formula is known for those numbers, they are equal to the number of certain semi-standard Young tableaux, which is their ``combinatorial interpretation''.

Looking at the analogous structure constants for $S_n$, the Kronecker coefficients $g(\la,\mu,\nu)$ give the multiplicities of $S_n$--irreducible representations $\mathbb{S}_\la$ in the tensor product of two others $\mathbb{S}_\mu \otimes \mathbb{S}_\nu$. Yet, despite their innocuous definition about 85 years ago, mimicking the Littlewood-Richardson one, no formula nor positive combinatorial interpretation is known for them. Likewise, no positive combinatorial interpretation is known for the plethysm coefficients of $GL_N$. 

But what is a ``combinatorial interpretation''? It would be a set of easily defined combinatorial objects, whose cardinality gives our desired structure constants. Yet, what do we consider combinatorial objects? For the most part, we know them once we see them, just like the LR tableaux. But could it be that no such nice objects exist, and, if so, could we prove that formally and save ourselves the trouble of searching for black cats in dark rooms, an endeavor particularly difficult when there is no cat. 

This is when Computational Complexity Theory comes into play and provides the formal framework we can settle these questions in. In its classical origins, Computational Complexity Theory classifies computational problems according to usage of resources (time and/or space) needed to obtain an answer. In our context, the resource is time as measured by the number of elementary steps needed to be performed by an algorithm solving that problem. Problems are thus divided into computational complexity classes depending on how fast they can be solved. For decision problems, that is when we are looking for a Yes/No answer, the main classes are $\P$, of problems solvable in polynomially (in the input size) many steps,  and $\NP$ is the class of problems, for which if the answer is Yes, then it could be verified in polynomial time. We have that $\P \subset \NP$ and the $\P$ vs $\NP$ Millennium problem asks whether $\P\neq \NP$, which is the widely believed hypothesis. The class $\NP$ is characterized by its complete problems like 3SAT or HAMCYCLE, which asks, given a graph $G$, whether it has a Hamiltonian cycle. If the graph has such a cycle, then one can specify it by its sequence of vertices, and verify in linear time that there is an edge between two consecutive ones.

For the counting problems, where the answer should be a nonnegative integer, the corresponding classes are $\FP$ and $\SP$. The class $\SP$ can also be defined as counting exponentially many polynomially computable discrete objects, and a $\SP$ formula  is a naturally nonnegative [exponentially large] sum of counting objects computable in polynomial time.

When there are ``nice'' formulas, like the hook-length formula, or the determinantal formulas for skew standard Young tableaux, the  corresponding problem (e.g. to compute $f^\la$) is in $\FP$. When we have a ``combinatorial interpretation'', the corresponding problem is in $\SP$, see~\cite{Pak22,Pan22}. Thus, to show that a ``reasonable combinatorial interpretation'' does not exist, we may want to show that the given problem is not in $\SP$ under some widely accepted assumptions, e.g. $\P \neq \NP$ or that the polynomial hierarchy $\PH$ does not collapse. 

In our quest to compute the Kronecker or plethysm coefficients, we ask whether the corresponding problem is in $\SP$. As a proof of concept we show that a similar problem, the computation of the square of an $S_n$ character, is not in $\SP$, given that the polynomial hierarchy does not collapse to second level,~\cite{IPP22}. 

\subsection{Geometric Complexity Theory}

In the opposite direction, Algebraic Combinatorics is used in Geometric Complexity Theory, a program aimed at finding computational lower bounds and distinguishing classes using Algebraic Geometry and Representation Theory. 

In his landmark paper \cite{V1} from 1979, Valiant defined algebraic complexity classes for computing polynomials in formal variables. Later these classes were denoted by $\VP$ and $\VNP$, and represented the algebraic analogues of the original $\P$ and $\NP$ classes.\footnote{As we shall see later, there is a fine distinction of the algebraic versions of $\P$: $\VPs, \VBP, \VP$, but it is often ignored. Formally, the relevant class is $\VBP$.} The flagship problem in arithmetic complexity theory is to show that $\VP \neq \VNP$ and is closely related to $\P \neq \NP$, see~\cite{Bur0}. As with $\P$~vs~$\NP$, the general strategy is to identify complete problems for $\VNP$, i.e.\ complete polynomials, and show they do not belong to $\VP$. Valiant identified such $\VNP$-complete polynomials, most notably the permanent of a $n\times n$ variable matrix. At the same time he showed that the determinant polynomial is $\VP$--universal, i.e.\ every polynomial from $\VP$ can be computed as a polynomially sized determinant of a matrix whose entries are affine linear forms in the original variables. This sets the general strategy of distinguishing $\VP$ from $\VNP$ by showing that the permanent is not a determinant of some poly-sized matrix.

Geometric Complexity Theory aims to distinguish such algebraic complexity classes via the algebro-geometric properties of the corresponding complete/universal polynomials.
In two landmark papers \cite{MS1, MS2} Mulmuley and Sohoni suggested distinguishing these polynomials by studying the algebraic varieties arising from the group action corresponding to all the linear transformations. In particular, to distinguish polynomials, one can consider the representation theoretic structure of these varieties ['s coordinate rings] and find some irreducible representations appearing with different multiplicities in the two. Because of the many symmetries of the polynomials and the equivariant action of $GL_N$, usually such multiplicities can  be naturally expressed via the fundamental structure constants Kronecker, Littlewood--Richardson, plethysms etc.\ from $\S$\ref{s:mult}, and the methods to study them revolve around the combinatorics of Young Tableaux and generalizations.

Since such multiplicities are even harder than the Kronecker and plethysm coefficient, a simpler approach would have been to study just the occurrence of irreducible representations rather than the value of the multiplicity. If an irreducible $GL$ module appears in the [coordinate ring orbit closure of the] permanent of an $m\times m$ matrix, but not for the determinant of an $n\times n$ matrix, that would imply a lower bound, namely the $\per_m$ cannot be equal to ${\det}_n$ (of affine linear forms). If this happens for $n>poly(m)$ (i.e. bigger than any polynomial in $m$), then $\VP \neq \VNP$. Such irreducible representations are called occurrence obstructions, and unfortunately do not exist~\cite{BIP} for this model. 

Thus we have to compare the actual multiplicities or explore other models besides permanent versus determinant. Understanding their growth starts with finding bounds and later asymptotics for Kronecker and plethysm coefficients, see~\cite{Pan22} for further discussions.

\subsection{Paper structure}

In Section~\ref{s:ac} we will define the basic objects in Algebraic Combinatorics and Representation Theory and recall important facts on SYTs and symmetric functions. In Section~\ref{s:mult} define the culprit structure constants Kronecker and plethysm coefficients and recall some of the major open problems. In Section~\ref{s:cct} we will discuss Computational Complexity Theory from the point of view of a mathematician. In Section~\ref{s:cc_in_ac} we will discuss how Computational Complexity can be applied in Algebraic Combinatorics, stating various hardness and completeness results and conjectures on Kostka, LR, Kronecker coefficients and the characters of the symmetric group. In Section~\ref{s:ac_in_gct} we will discuss Geometric Complexity Theory in more detail, explain the connection with Algebraic Combinatorics and some of the recent advances in the area. The text will aim to be as self-contained as possible. 

For other open problems on structure constants, in particular positivity and asymptotics, see~\cite{Pan22}. 

\smallskip

\noindent{\bf Disclaimer.} The current paper is a detailed transcript of the author's Current Developments in Mathematics talks (April 2023). This work does not attempt to be a broad survey on the topic, and is naturally concentrated on the author's viewpoint and work. 

\bigskip

\noindent{\bf Acknowledgments.} The author is grateful to Christian Ikenmeyer and Igor Pak for the years of fruitful collaborations on the subject addressed here. Many thanks also to Sara Billey,  Allen Knutson,   Alex Yong for many useful questions and discussions   on these topics. The author has been partially supported by the NSF.

\section{Algebraic Combinatorics}\label{s:ac}

Here we will describe the basic objects and facts from Algebraic Combinatorics which will be used later. For further details on the combinatorial sides see~\cite{S1,Mac} and for the representation theoretic aspects see~\cite{Sag,Ful}.

\subsection{Partitions and Tableaux}\label{ss:bakcground}

\ytableausetup{boxsize=2ex}
\defn{Integer partitions} $\la$ of $n$, denoted $\la \vdash n$, are sequences of nonnegative integers $\la=(\la_1,\la_2,\ldots,\la_k)$, such that $\la_1 \geq \la_2 \geq \cdots \geq 0$ and $\la_1+\cdots +\la_k=n$. We denote by $\ell(\la) = \max\{i: \la_i >0\}$ the \defn{length} of the partition, which is the number of nonzero parts, and by $|\la| = \la_1+\la_2+\cdots +\la_\ell$ its \defn{size}. A partition can be represented as a \defn{Young diagram}, which is a left-justified array of squares, such that row $i$ (indexing top to bottom) has exactly $\la_i$ squares. For example $\la=(4,3,1) \vdash 8$ has $\ell(\la)=3$ and its Young diagram is \ydiagram{4,3,1}. Here we will denote by $[\la]$ the Young diagram and think of it as a set of squares with coordinates $(i,j)$ with $(1,1)$ being the topmost leftmost box, so the box at $(2,3)$ is $\ydiagram{4,3,1}*[*(blue)]{3+0,2+1}$. We denote by $(1^k) = (\underbrace{1,\ldots,1}_k)$ the single column partition of $k$ boxes,  call \emph{one-row} partition the partitions with only one nonzero part, \emph{two-row} partitions of the form $(n-k,k)$, and \emph{hooks} partitions of the kind $(n-k,1^k)$. We denote by $(a^k) = (\underbrace{a,\ldots,a}_k)$ the \emph{rectangular} partition whose Young diagram is a $k \times a$ rectangle.
The transpose or \defn{conjugate} partition of $\la$ is denoted $\la'$ and is the one whose shape is obtained by transposing along the main diagonal $[\la]$, e.g. for $\la=(4,3,1)$ we have $\la'=(3,2,2,1)$. The \defn{skew} partition $\la/\mu$ is obtained by removing the squares occupied by $\mu$ from $\la$, so for example the Young diagram of $(5,4,3,1) /(2,1)$ is $\ydiagram{2+3,1+3,3,1}$.  

The set of partitions of $n$ will be denoted by $\mathcal{P}(n)$ and its cardinality by $p(n)$. While there is no closed form formula for $p(n)$, there is a nice generating function
$$ \sum_{n=0}^\infty p(n) t^n = \prod_{i=1}^\infty \frac{1}{1-t^i}.$$ 

A \defn{standard Young tableaux} (SYT) of \defn{shape} $\la\vdash n$  is a bijection $T:[\la] \xrightarrow{\sim} \{1,\ldots,n\}$, such that $T(i,j) <T(i+1,j)$ and $T(i,j) < T(i,j+1)$. For example, the SYTs of shape $(2,2,1)$ are
$$\ytableaushort{12,34,5} \qquad \ytableaushort{12,35,4} \qquad \ytableaushort{13,24,5} \qquad \ytableaushort{13,25,4} \qquad \ytableaushort{14,25,3}\,.$$

The \defn{hook} of a box $(i,j)$ in $[\la]$ is the collection of squares $\{ (i,j), (i+1,j),\ldots, (\la'_j,j), (i,j+1),\ldots,(i,\la_i)\}$ below $(i,j)$ in the same column, or to the right in the same row. For example, for box $(2,2)$ in $(5,4,3,3)$ the hook is $\ydiagram{5,4,3,3}*[*(cyan)]{0+0,1+3,1+1,1+1}$. The hook-length $h_{i,j}$ is the number of boxes in the hook of $(i,j)$.

Let $f^\la$ be the number of standard Young tableaux of shape $\la \vdash n$. Then the \emph{hook-length formula} (HLF) of Frame-Robinson-Thrall~\cite{FRT} gives
\begin{equation}\label{eq:hlf}
f^\la = \frac{n!}{\prod_{(i,j) \in [\la]}  h_{i,j} }.
\end{equation}

The following remarkable identity
\begin{equation}\label{eq:rsk}
\sum_\la (f^{\la})^2 = n!
\end{equation}
gives rise to an even more remarkable bijection between pairs of same shape SYTs and permutations, known as RSK for Robinson-Schensted-Knuth.
For example 
 $$\left( \ytableaushort{124,3}, \ytableaushort{123,4} \right) \longleftrightarrow  4123\ .$$

The \defn{semi-standard Young tableaux} (SSYT) of shape $\la$ and content $\al$ are maps $T:[\la] \to \mathbb{N}$, such that $|T^{-1}(i)|=\al_i$, i.e. $\al_i$ many entries are equal to $i$, and the entries increase weakly along rows and strictly down columns, i.e. $T(i,j) \leq T(i,j+1)$ and $T(i,j) < T(i+1,j)$. We denote the set of such tableaux by $\SSYT(\la;\al)$.

For example the SSYTs of shape $\la=(3,3,1)$ and type $\al=(2,2,2,1)$ are

$$\ytableaushort{112,233,4} \qquad \ytableaushort{112,234,3} \qquad \ytableaushort{113,224,3}\ . $$

Completely analogously, we can define the \defn{skew SSYT} of shape $\la/\mu$ as the fillings of $[\la/\mu]$ with integers weakly increasing along rows and strictly down columns, e.g. 
$\ytableaushort{ \none\none 23,\none 144,23}$
is a skew SSYT of shape $(4,4,2)/(2,1)$ of type $(1,2,2,2)$. 

\subsection{Symmetric functions}

Let $\Lambda[\bx]$ be the ring of \defn{symmetric functions} $f(x_1,x_2,\ldots)$,
where the symmetry means that $f(\bx)=f(\bx_\sigma)$ for any permutation $\sigma$
of the variables, and $f$ is a formal power series. When all but  finitely many variables are 0 then $f$ becomes a symmetric polynomial. 

The ring $\Lambda$ is graded by the total degree, and its component of degree $n$ has dimension $p(n)$ as a $\mathbb{C}$-vector spaces. There are several useful bases for $\Lambda$ -- the monomial, the elementary, power sum, (complete) homogeneous, and Schur functions.

 The
\emph{monomial symmetric functions} $m_\la (x_1,\ldots,x_k)$ are defined as the sum of all distinct
monomials of the form \ts $x_{\si(1)}^{\la_1}\cdots x_{\si(k)}^{\la_k}$, where $\si \in S_k$.
For example, \ts $m_{311}(x_1,\ldots,x_n) = x_1^3x_2x_3 + x_1^3x_2x_4 + \ldots$, and each monomial appears with coefficient 0 or 1. 

Let $p_k := m_{(k)} = x_1^k + x_2^k + \ldots$. The \emph{power sum symmetric functions}
are then defined as \. $p_\la := p_{\la_1}\. p_{\la_2} \cdots$

The {\em elementary symmetric functions} $\{e_\la\}$ are defined  as follows
$$e_k := \sum_{i_1<i_2<\cdots <i_k} x_{i_1} x_{i_2} \cdots x_{i_k} , \qquad \text{ and }\qquad e_\la := e_{\la_1}e_{\la_2}\cdots.$$
For example 
$$e_{2,1}(x_1,x_2,x_3) = (x_1x_2+x_1x_3 +x_2x_3)(x_1+x_2+x_3) = m_{2,1}(x_1,x_2,x_3) + 3m_{1,1,1}(x_1,x_2,x_3).$$

The  {\em homogeneous symmetric functions} $h_\la$ are given by 
$$h_k := \sum_{i_1\leq i_2\leq \cdots \leq i_k} x_{i_1} x_{i_2} \cdots x_{i_k} , \qquad \text{ and }\qquad h_\la := h_{\la_1}h_{\la_2}\cdots.$$

 The \defn{Schur functions} can be defined as the  generating functions of SSYTs, where $\bx^\al:=x_1^{\al_1}\cdots x_k^{\al_k}$, namely
 $$s_\la = \sum_{\al} \sum_{T \in \SSYT(\la,\al)} \bx^\al,$$
 where $\al$ goes over all weak compositions of $n$. 
For example, $s_{1^k} = e_k$, $s_{k} = h_k$ and 
$$s_{(2,1)}(x_1,x_2,x_3) = x_1^2x_2 + x_1 x_2^2 +x_1^2x_3 + x_1x_3^2 + x_2^2x_3 +x_2x_3^2 + 2x_1x_2x_3.$$
We can also define the \emph{skew Schur functions} $s_{\la/\mu}$ as the analogous generating function for $\SSYT(\la/\mu)$.

They can also be defined and computed via the \emph{Weyl determinantal formula} 
 
$$
s_\la(x_1,\ldots,x_k) \, := \, \det\Bigl( x_i^{\lambda_j+k-j}\Bigr)_{i,j=1}^k
\prod_{1\le i<j\le n} \frac{1}{x_i-x_j} 
$$
or the \emph{Jacobi-Trudi identity}
$$
s_{\la/\mu} = \det [ h_{\la_i  -i - \mu_j+j} ]_{i,j=1}^{\ell(\la)}.$$

The ring $\Lambda$ has an inner product $\langle \cdot , \cdot \rangle$, where the Schur functions form an orthonormal basis and the power sums are orthogonal. Namely

$$\langle s_\la, s_\mu \rangle = \delta_{\la,\mu} \qquad \langle p_\la, p_\mu \rangle = z_\la \delta_{\la,\mu}.$$
Additionally $\langle h_\la, m_\mu \rangle = \delta_{\la,\mu}$, where $\delta_{\la,\mu}=0$ if $\la \neq \mu$ and $1$ if $\la=\mu$ and $z_\la = \frac{n!}{\prod_i i^{m_i}m_i!}$ when $\la = (1^{m_1}2^{m_2}\cdots)$, i.e. there are $m_i$ parts equal to $i$.  

The involution $\omega$ is defined as $\omega(e_\la) = h_\la$, $\omega^2=id$ and we have that $\omega(s_\la)=s_{\la'}$.

The Schur functions posses beautiful combinatorial properties, for example they satisfy the \emph{Cauchy identity}
$$\sum_{\la} s_\la(\bx) s_\la(\by) = \prod_{i,j} \frac{1}{1-x_iy_j} ,$$
which can also be proven via RSK.

\subsection{Representations of $S_n$ and $GL_N$}
A  \defn{group representation} $\rho$ of a group $G$ is a group homomorphism   $\rho:G \to GL(V)$, which can also be interpreted as an action of the group on a vector space $V$. We often refer to the vector space $V$ as the representation. An irreducible representation is such a vector space $V$ which has no nontrivial invariant subspaces. If $G$ is finite or reductive and the underlying field is $\mathbb{C}$ then every representation can be uniquely decomposed as direct sum of irreducible representations. Such decompositions can be easily studied via the characters, $\chi^{\rho}(g):= trace( \rho(g) )$, which are central functions, i.e. constant on conjugacy classes.

The irreducible representations of the \emph{symmetric group} $S_n$ are the \defn{Specht modules} $\mathbb{S}_\la$ and are indexed by partitions $\la \vdash n$. Using row and column symmetrizers in the group algebra $\mathbb{C}[S_n]$ one can construct the irreducible modules as certain formal sums over tableaux of shape $\la$. Each such element has a unique minimal tableau which is an SYT, and so a basis for $\mathbb{S}_\la$ can be given by the SYTs. In particular
$$\dim \mathbb{S}_\la = f^\la.$$

We have that $\mathbb{S}_{(n)}$ is the trivial representation assigning to every $w$ the value $1$ and $\mathbb{S}_{1^n}$ is the sign representation. 

The \defn{character $\chi^\la(w)$} of $\mathbb{S}_{\la}$ can be computed via the \emph{Murnaghan-Nakayama} rule. Let $w$ have type $\al$, i.e. it decomposes into cycles of lengths $\al_1,\al_2,\ldots,\al_k$. Then
$$\chi^\la(w)=\chi^\la(\al) = \sum_{T \in MN(\la;\al)} (-1)^{ht(T)},$$
where $MN$ is the set of rim-hook tableaux of shape $\la$ and type $\al$, so that the entries are weakly increasing along rows and down columns, and all entries equal to $i$ form a rim-hook shape (i.e. connected, no $2\times 2$ boxes) of length $\al_i$. The height of each rim-hook is one less than the number of rows it spans, and $ht(T)$ is the sum of all these heights. For example,
$$\ytableaushort{112333,12234,22334}$$
is a Murnaghan-Nakayama tableau of shape $(6,5,5)$, type $(3,5,6,2)$ and has height $ht(T) = 1 + 2 + 2 +1=6$.

As we shall see, the characters are also the transition matrices between the $\{s_\la\}$ and $\{p_\la\}$ bases of $\Lambda$. 

The irreducible polynomial representations of  $GL_N(\mathbb{C})$ are the \emph{Weyl modules} $V_\la$ and are indexed by all partitions with $\ell(\la) \leq N$. Their characters are exactly the Schur functions $s_\la(x_1,\ldots,x_N)$, where $x_1,\ldots,x_N$ are the eigenvalues of $g \in GL_N(\mathbb{C})$. The dimension is just $s_\la(1^N)$ and can be evaluated as a product via the \emph{hook-content formula}
$$s_\la(1^N) = \prod_{(i,j) \in [\la] } \frac{N + j-i}{h_{i,j}}.$$

\section{Multiplicities and structure constants}\label{s:mult}

\subsection{Transition coefficients}\label{ss:transition_coeff}

As the various symmetric function families form bases in $\Lambda$, it is natural to describe the transition matrices between them. The coefficients involved have significance beyond that.

We have that 
$$h_\la = \sum_{\mu} CT(\la,\mu) m_\mu,$$
where $CT(\la,\mu)$ is the number of \emph{contingency arrays} $A$ with marginals $\la, \mu$, namely $A \in \mathbb{N}_0^{\ell(\la)\times \ell(\mu)}$ and $\sum_j A_{ij} = \la_i$, $\sum_i A_{ij} = \mu_j$. 

Similarly, 
$$e_\la = \sum_{\mu}  CT_0(\la,\mu) m_{\mu},$$
where $CT_0(\la,\mu)$ is the number of $0-1$ contingency arrays $A$ with marginals $\la,\mu$, i.e. $A \in \{0,1\}^{\ell(\la)\times \ell(\mu)}$.

We have that 
$$p_\la = \sum_{\mu} P(\la,\mu) m_\mu,$$
where for any two integer vectors $\mathbf{a},\mathbf{b}$, we set
$$P( \mathbf{a},\mathbf{b}) := \# \{  (B_1,B_2,\ldots,B_k): B_1 \sqcup B_2 \sqcup \ldots \sqcup B_k =[m], \sum_{i  \in B_j}  a_i = b_j \text{ for all }j=1,\ldots,k\} $$
be the number of \defn{ordered set partitions} of items $\mathbf{a}$ into bins of sizes $\mathbf{b}$. 

The \defn{Kostka numbers} $K_{\la\mu}$, $\la,\mu\vdash n$ are defined by
$$
h_\la \, = \, \sum_{\mu \vdash n} \, K_{\la \mu} \, s_\la$$
and by orthogonality also as
$$s_\la \, = \, \sum_{\mu\vdash n} \, K_{\la \mu} \, m_\mu\,.
$$
By definition we have that $L_{\la,\mu}=|\SSYT(\la,\mu)|$, i.e. the number of SSYTs of shape $\la$ and type $\mu$.

Finally, the symmetric group characters appear as 
$$p_\al = \sum_\la \chi^\la(\al)s_\la$$
or equivalently, as
$$s_\la = \sum_\al \chi^\la(\al) z_\al^{-1} p_\al.$$

\subsection{Tensor products and structure constants}\label{ss:str_const}

Once the irreducible representations have been sufficiently understood, it is natural to consider what other representations can be formed by them and how such representations decompose into irreducibles. 
Such problems are often studied in quantum mechanics under the name \emph{Clebsh-Gordon coefficients}. 

In the case of $GL_N(\mathbb{C})$ these coefficients are
the  \defn{Littlewood--Richardson} \defn{coefficients} (LR) $c^{\la}_{\mu\nu}$ defined as the multiplicity of $V_\la$ in $V_\mu \otimes V_\nu$, so
$$V_\mu \otimes V_\nu = \bigoplus_\la V_\la ^{\oplus c^{\la}_{\mu\nu}}.$$
Via their characters, they can be equivalently defined as

$$s_\mu s_\nu = \sum_\la c^{\la}_{\mu\nu} s_{\la}, \quad \text{ and } \quad s_{\la/\mu} = \sum_\nu c^{\la}_{\mu\nu}s_\nu.$$

While no nice product formula exists for their computation, they have a \emph{combinatorial interpretation}, the so called \emph{Littlewood-Richardson rule}. This rule was first stated by Littlewood and Richardson in 1934~\cite{LR34}, survived through several incomplete proofs, until formally proven using the new technology listed here by Sch\"utzenberger and, separately, Thomas in the 1970s.  

\begin{thm}[Littlewood-Richardson rule]
The Littlewood-Richardson coefficient $c^{\la}_{\mu\nu}$ is equal to the number of skew SSYT $T$ of shape $\la/\mu$ and type $\nu$, whose reading word is a ballot sequence.
\end{thm} 
The \defn{reading word} of a tableau $T$ is formed by reading its entries row by row from top to bottom, such that the rows are read from the back. For example the reading word of $\ytableaushort{\none \none 112,\none233,14}$\, is $21133241$. A sequence $a_1a_2\ldots$ is a \defn{ballot sequence} if for every $i$ and every $k$ we have $\#\{j: a_j=i, j\leq k\} \geq \#\{j: a_j=i+1, j\leq k\}$, i.e. there are weakly more $i$'s thant $(i+1)$'s in every initial segment $a_1\ldots a_k$ of the sequence. So the above example is not a ballot sequence because for $k=1$ we have more $2$s than $1$s.

\begin{ex} We have $c^{(6,4,3)}_{(3,1),(4,3,2)}=2$ as there are two LR tableaux (SSYT of shape $(6,4,3)/(3,1)$ and type $(4,3,2)$ whose reading words are ballot sequences)
\ytableaushort{\none\none\none 111,\none 122,233}\, with reading word $111221332$ and \ytableaushort{\none\none\none 111,\none 222,133}\,   with reading word 111222331. 
\end{ex}

It is easy to see the similarity with the Kostka numbers, and indeed, Kostka is a special case of LR via the following
\begin{equation}\label{eq:kostka-lr}
K_{\la,\mu} = c^{\theta}_{\la_1^{\ell-1}, \tau},
\end{equation}
where $\ell = \ell(\mu)$, $\theta = (\la_1^{\ell-1} + \eta,\la)$ with $\eta = (\mu_1*(\ell-1), \mu_1*(\ell-2),\ldots,\mu_1)$ and $\tau = \eta + \mu$. 

\begin{ex}
Let $\la=(3,2)$ and $\mu=(2,2,1)$, then the suggested $LR$ tableaux by the above formula would be 
\ytableaushort{\none\none\none 1111,\none\none\none 22,113,22}\, and \ytableaushort{\none\none\none 1111,\none\none\none 22,112,23}.
As we see the regular SSYTs of shape $\la$ and type $\mu$ emerge in the bottom. The top parts are forced and their reading words have many more 1s than 2s, more 2s than 3s etc so that they overwhelm the ballot and the ballot condition is trivially satisfied by any SSYT in the bottom part. 
\end{ex}

The \defn{Kronecker coefficients} $g(\la,\mu,\nu)$ of the symmetric group are the corresponding structure constants for the ring of $S_n$- irreducibles. Namely, $S_n$ acts diagonally on the tensor product of two Specht modules and the corresponding module factors into irreducibles with multiplicities $g(\la,\mu,\nu)$
$$\mathbb{S}_\la \otimes \mathbb{S}_\mu = \bigoplus_\nu \mathbb{S}_\nu^{\oplus g(\la,\mu,\nu)}, \text{ i.e.} \quad \chi^\la \chi^\mu = \sum_\nu g(\la,\mu,\nu) \chi^\nu
$$
In terms of characters we can write them as 
\begin{equation}\label{eq:char_kron}
g(\la,\mu,\nu) = \langle \chi^\la \chi^\mu,\chi^\nu\rangle = \frac{1}{n!} \sum_{w\in S_n} \chi^\la(w) \chi^\mu(w) \chi^\nu(w).
\end{equation}
The last formula shows that they are symmetric upon interchanging the underlying partitions $g(\la,\mu,\nu)= g(\mu,\nu,\la) = \cdots$ which motivates us to use such symmetric notation. 

The \defn{Kronecker product} $*$ on $\Lambda$ is defined on the Schur basis by 
$$s_\la * s_\mu = \sum_\nu g(\la,\mu,\nu) s_\nu,$$
and extended by linearity. 

The Kronecker coefficients were defined by Murnaghan in 1938~\cite{Mur38}, who was inspired by the Littlewood-Richardson story. In fact, he showed that
\begin{thm}[Murnaghan]\label{thm:mur}
For every $\la,\mu,\nu$, such that $|\la|=|\mu|+|\nu|$ we have that 
$$c^\la_{\mu\nu} = g( (n - |\la|,\la), (n-|\mu|,\mu), (n-|\nu|,\nu)),$$
for sufficiently large $n$.
\end{thm}
In particular, one can see that $n=2|\la|+1$ would work. Note that this also implies that these Kronecker coefficients stabilize as $n$ increases, which is also true in further generality. 

Thanks to the Schur-Weyl duality, they can also be interpreted via Schur functions as
$$s_\la[\bx \by] = \sum_{\mu,\nu} g(\la,\mu,\nu) s_\mu(\bx) s_\nu(\by),$$
where $\bx\by = (x_1y_1,x_1y_2,\ldots,x_2y_1,\ldots)$ is the vector of all pairwise products. In terms of $GL$ representations they give us the dimension of the invariant space
$$g(\la,\mu,\nu) = \dim ( V_\mu \otimes V_\nu \otimes V^*_\la)^{GL_{N}\times GL_M},$$
where $V_\mu$ is considered a $GL_N$ module, and $V_\nu$ a $GL_M$ module.
Using this interpretation for them as dimensions of highest weight spaces, see~\cite{CHM}, one can show the following
\begin{thm}[Semigroup property]\label{thm:semigroup}
Let $(\al^1,\be^1,\ga^1)$ and $(\al^2,\be^2,\ga^2)$ be two partition triples, such that $|\al^1|=|\be^1|=|\ga^1|$ and $|\al^2|=|\be^2|=|\ga^2|$. Suppose that $g(\al^1,\be^1,\ga^1)>0$ and $g(\al^2,\be^2,\ga^2)>0$. Then
$$g(\al^1+\al^2,\be^1+\be^2,\ga^1+\ga^2) \geq \max\{ g(\al^1,\be^1,\ga^1), g(\al^2,\be^2,\ga^2)\}.$$
\end{thm}

Other simple properties we can see using the original $S_n$ characters are that 
$$g(\la,\mu,\nu) = g(\la',\mu',\nu)$$
since $\mathbb{S}_\la = \mathbb{S}_{\la'} \otimes \mathbb{S}_{1^n}$, where $\chi^{1^n}(w) =\sgn(w)$ is simply the sign representation. Similarly, we have $g(\la,\mu,(n))=\delta_{\la,\mu}$ for all $\la$ and $g(\la,\mu,1^n)=\delta_{\la,\mu'}$. 

\begin{ex}
By the above observation we have that  
$$h_k[\bx\by] =s_k[\bx\by] = \sum_{\la \vdash k} s_\la(\bx) s_\la(\by).$$
Using the Jacobi-Trudi identity we can write
\begin{align*}
s_{2,1}[\bx\by] = h_2[\bx\by]h_1[\bx\by] - h_3[\bx\by] = \left(s_2(\bx)s_2(\by) +s_{1,1}(\bx)s_{1,1}(\by)\right)s_1(\bx)s_1(\by) \\
- s_3(\bx)s_3(\by) -s_{2,1}(\bx)s_{2,1}(\by) - s_{1,1,1}(\bx)s_{1,1,1}(\by)\\
= s_{2,1}(\bx)s_{2,1}(\by) + s_{2,1}(\bx)s_3(\by) +s_3(\bx)s_{2,}(\by) +s_{1,1,1}(\bx)s_{2,1}(\by) + s_{2,1}(\bx)s_{1,1,1}(\by).
\end{align*}
So we see that $g((2,1),(2,1),(2,1))=1$. 
\end{ex}

The \defn{plethysm coefficients} $a^\la_{\mu,\nu}$ are multiplicities of an irreducible $\GL$ representation in the composition of two $\GL$ representations. Namely, let $\rho^\mu: GL_N \to GL_M$ be one irreducible, and $\rho^\nu: GL_M \to GL_K$ be another. Then $\rho^\nu \circ \rho^\mu: GL_N \to GL_K$ is another representation of $GL_N$ which has a character $s_\nu[s_\mu]$, which decomposes into irreducibles as
$$s_\nu[s_\mu] = \sum_\la a^\la_{\mu,\nu}s_\la.$$
Here the notation $f[g]$ is the evaluation of $f$ over the monomials of $g$ as variables, namely if $g = \bx^{\al^1} + \bx^{\al^2} +\cdots$, then $f[g] = f(\bx^{\al^1},\bx^{\al^2},\ldots)$. 

\begin{ex}
We have that
\begin{align*}
s_{(2)}[s_{(1^2)}] = h_3[e_2] = h_2(x_1x_2,x_1x_3,\ldots)=x_1^2x_2^2 + x_1^2x_2x_3 + 3x_1x_2x_3x_4+\cdots \\
= s_{2,2}(x_1,x_2,x_3,\ldots) + s_{1,1,1,1}(x_1,x_2,x_3,\ldots),  
\end{align*}
so $a^{(2,2)}_{(2),(1,1)} =1$ and $a^{(3,1)}_{(2),(1,1)}=0$.
\end{ex}

We will be particularly interested when $\mu = (d)$ or $(1^d)$ which are the $d$th symmetric power $Sym^d$ and the $d$th wedge power $\Lambda^d$, and $\nu = (n)$. We denote this plethysm coefficient by $a_\la(d[n]) := a^\la_{(d),(n)}$ and
$$h_d[h_n] = \sum_\la a_\la(d[n]) s_\la.$$

The following can easily be derived using similar methods, see~\cite{PP14}.
\begin{prop}
We have that $g(\la,n^d, n^d) = a_\la(d[n])=p_{\la_2}(n,d) - p_{\la_2-1}(n,d)$ for $\la \vdash nd$, such that $\ell(\la)\leq 2$. Here $p_r(a,b) = \# \{ \mu \vdash r: \mu_1\leq a, \ell(\mu)\leq b\}$ are the partitions of $r$ which fit inside a rectangle. 
\end{prop}

In particular, these are the coefficients in the $q$-binomials
 $$\sum_r p_r(a,b) q^r = \binom{a+b}{a}_q := \prod_{i=1}^a \frac{ (1-q^{i+b})}{1-q^i}.$$
 
 As a curious application, these identities were used in~\cite{PP13,PPq} to prove the strict version of Sylvester's unimodality theorem and find bounds on the coefficients of the $q$-binomials. Later in \cite{MPPe}, using tilted random geometric variables, we found tight asymptotics for the differences of $p_r(a,b)$ and hence obtained tight asymptotics for this family of Kronecker coefficients.

\section{Computational Complexity Theory}\label{s:cct}

Here we will define, in broad and  not fully precise terms, the necessary computational complexity classes and models of computation. For background on the subject we refer to~\cite{Aa, Bur00, Wig}. 

\subsection{Decision and counting problems}

Computational problems can be classified depending on how much of a given resource (time or memory) is needed to solve it via an algorithm, i.e. produce the answer for any given input of certain size. Depending on the model of computation used (e.g. Turing machines, boolean circuits, quantum computers etc) the answers could vary. Here we will only focus on classical computers and will consider complexity depending on the time an algorithm takes, which is essentially equivalent to the number of elementary steps an algorithm performs. 

Let $I$ denote the input of the problem and let $|I|=n$ be its size as the number of bits it takes to write down in the computer. Depending on the encoding of the problem the size can vary, and then the ``speed'' of the algorithm as a function of the size will change. There are two main ways to present an input: \defn{binary} versus \defn{unary}. If the input is an integer $N$, then in binary it would have size about $\lceil \log_2(N) \rceil$, for example when $N=2023$, in binary it is $11111100111$ and the input size is $11$. In unary, we will write $N$ as $\underbrace{111\ldots1}_N$ and in our case take up $N=2023$ bits. 
As we shall see soon, the encoding matters significantly on how fast the algorithms are as functions of the input size.  From complexity standpoint, encoding in binary or in any other base $b>1$, does not make a difference as the input size is just rescaled $\log_b N = \log_b(2) \log_2 N $.

A \defn{decision problem}, often referred to as \emph{language}, is a problem, whose answer should be Yes/No. For example, in the problem $\textsc{PRIMES}$ we are given input $N$ and have to output Yes if $N$ is a prime number. The \defn{complexity class $\P$} consists of the decision problems, such that the answer can be obtained in \defn{polynomial time}, that is $O(n^k)$ for some fixed $k$ (fixed for the given problem, but independent of the input). Thanks to the~\cite{AKS} breakthrough result, we now know  that $\textsc{PRIMES} \in \P$.

The \defn{complexity class $\NP$} consists of decision problems, such that if the answer is Yes then it can be \emph{verified} in polynomial time, i.e. they have a \emph{poly-time witness}.  The problem is phrased as ``given input $I$, is the set $C(I)$ nonempty''. It is in $\NP$  iff whenever $C(I) \neq \emptyset$, then there would be an element $X(I) \in C(I)$, such that we can check whether $X(I) \in C(I)$ in $O(n^k)$ time for some fixed $k$. For example, in $\textsf{HAMCYCLE}$, the input is a graph $G=(V,E)$ (encoded as its adjacency matrix, so the input size is $O(|V|^2)$), and the question is ``does $G$ have a Hamiltonian cycle''. In this case $C(G)$ would be the set of all Hamiltonian cycles in $G$ (encoded as permutations of the vertices), and given one such cycle $X(G)=v_1\ldots v_m$  we can check in $O(m)$ time whether it is indeed a Hamiltonian cycle by checking whether $(v_i,v_{i+1}) \in E$ for all $i=1,\ldots,m$. 

We say that a problem is \defn{$\NP$-complete} if it is in $\NP$ and every other problem from $\NP$ can be reduced to it in polynomial time. A set of $\NP$-complete problems includes \textsf{HAMCYCLE, 3SAT, BINPACKING} etc. A problem is \defn{$\NP$-hard} if  every problem in $\NP$ is reducible to it in poly time. 

\begin{ex}
Here is an example when input size starts to matter. The problem \textsf{KNAPSACK} is as follows:

 Given an input $a_1,\ldots,a_m,b$ of $m+1$ integers, determine whether there is a subset $S \subset\{1,\ldots,m\}$, such that $\sum_{i \in S} a_i=b$. If the input integers $a_i$ are encoded in binary then the problem is $\NP$-complete. However, if they are encoded in unary then there is a dynamic programming algorithm that would output the answer in polynomial time. It is said that such problems can be solved in {\em pseudopolynomial time}. However the modern treatment would consider these problems as two different computational problems, one for each input encoding.
\end{ex}

We have that $\P \subset \NP$, but we are nowhere near showing that the containment is strict.

\begin{problem}[The $\P$ vs $\NP$ Millennium problem]
Show that $\P \neq \NP$.
\end{problem}

However, most researchers believe (and assume) that $\P \neq \NP$.

The \defn{class \coNP} consists of the decision problems, such that if the answer is No, then there exists a poly-time witness proving that. $X \in \coNP$ if and only if $\overline{X} \in \NP$. For example, $\overline{\textsf{HAMCYCLE}}$ would be the problem of deciding whether a graph does NOT have a Hamiltonian cycle. If the answer is no, then the graph has such a cycle and we can check it as above.  
The \defn{polynomial hierarchy} $\PH$ is a generalization of $\NP$ and $\coNP$ and is, informally, the set of all problems which can be solved by some number of augmentations by an oracle. Specifically, denote by $B^A$ the set of problems which can be solved by an algorithm from $B$ augmented with an ``oracle'' (another machine) from $A$. 
Then  we set $\Delta_0^\P:=\Sigma_0^\P:=\Pi_0^\P:=\P$ and recursively $\Delta^\P_{i+1} :=\P^{\Sigma_i^\P}$, $\Sigma_{i+1}^\P:=\NP^{\Sigma_i^\P}$ and $\Pi_{i+1}^\P:=\coNP^{\Sigma_i^\P}$. We set $\PH := \cup_i \left( \Sigma_i^\P \cup \Pi_i^\P \cup \Delta_i^\P\right)$. We have that $\Sigma_i \subset \Delta_{i+1} \subset \Sigma_{i+1}$ and $\Pi_i \subset \Delta_{i+1} \subset \Pi_{i+1}$, and it is yet another big open problem to prove the containments are strict. A widely believed hypothesis is that $\NP \neq \coNP$ and that $\PH$ does not collapse to any level (i.e. $\Sigma_i \neq \Sigma_{i+1}$ etc).

\smallskip

\defn{Counting problems} ask for the number of elements in $C(I)$ given input $I$. There are  two main complexity classes \FP and \SP, also believed to be different. \FP is the class of problems, such that $|C(I)|$ can be found in $O(n^k)$ time for some  fixed $k$. \defn{The class $\SP$} is the counting analogue of $\NP$ and can be defined as
$\SP$ is the class of functions $f: \{0,1\}^* \to \mathbb{N}$, such that there exists a polynomial $p$ and a verifier $V$ so that for an input $I$ we have 
\begin{align}\label{eq:sp_formula}
f(I) = |\{ y \in \{0,1\}^{p(|I|)}: V(I,y)=1\}|=\sum_{y=0}^{2^{p(|I|)}-1} V(I,y) .
\end{align}
 The verifier should be an algorithm running in polynomial time. That is, $\SP$ is the set of functions $f$ which can be expressed as 
\emph{exponentially large sums of  terms $V \in\{0,1\}$ which are determined in polynomial time}. 

\begin{ex}
\textsf{$\#$PERFECTMATCHINGS, $\#$HAMCYCLES, $\#$SETPARTITIONS} are all \SP-complete problems. In the case of \textsf{HAMCYCLE}, we have the input $I:=G$ a graph on $m$ vertices and $|I|=O(m^2)$ given by the edge pairs. Then $f$ counts the number of Hamiltonian cycles, so it can be computed by going  over all $m! = O(2^{m \log m})$ permutations $y:=v_{\sigma}$ of the vertices $\{v_1,\ldots,v_m\}$ and the verifier is $V(G,v_\sigma)=1$ iff $(v_{\si(i)},v_{\si(i+1)}) \in E(G)$ is an edge for every $i$. 
\end{ex}

The \defn{class $\GapP$} is the class of problems which are the difference of two $\SP$ functions, namely
$\GapP = \{ f-g: f,g \in \SP\}$. The class $\GapP_{\geq 0} = \GapP \cap \{f \geq 0\}$ is the class of $\GapP$ functions whose values are nonnegative. We define \defn{$\CeqP = [\GapP=0]$}, the  class of decision problems on whether two $\SP$ functions are equal. 

\medskip

The application of Computational Complexity theory in Combinatorics revolves around the following paradigm, see~\cite{Pak22} for detailed treatment. 
\smallskip

\begin{tabular}{p{3in}p{3in}}
\hline
 
{\em Counting and characterizing 
 combinatorial objects  given input data $I$ }  & {\em Solve: is $X \in C(I)$? or compute $|C(I)|$ } \smallskip \\
 \hline
 \smallskip 
 ``Nice formula'' & \smallskip The problem is in $\P, \FP$ \smallskip \\
 \hline
 \smallskip Positive combinatorial formula & \smallskip The problem is in $\NP$, $\SP$ \smallskip \\ 
 \hline
 \smallskip
No ``combinatorial interpretation'' & \smallskip The problem is not in $\SP$ \\
\hline
\end{tabular} 

\begin{rem}
The class $\SP$ is quite large. While it contains essentially all positive combinatorial formulas/interpretations we encounter in practice, it may actually be too large and other complexity classes like \ComCla{AC} could be more appropriate for certain problems.
\end{rem}

\begin{rem}
The above table does not address how the input is presented, but we can argue that there are {\em natural} encodings for the problems we will consider. Namely, we will see that if the inputs are in unary of input size $n$ then our problems of interest are in $\GapP$ and are nonnegative functions. This makes it natural to ask for the positive combinatorial formula to be in $\SP$. Moreover, the bounds on the sizes of our answers would be at most exponential in the input size $n$, i.e. $O(2^{p(n)})$ for some fixed degree polynomial $p$, which again suggests that a positive combinatorial formula is exactly of the kind~\eqref{eq:sp_formula}. Thus, problems like ``the number of sets of subsets of a given set'' are excluded from this consideration being ``too large'' for their input size. 
\end{rem}

\medskip

Besides the classical computational complexity, there is also \defn{quantum complexity}, informally defined by the minimal size of a quantum circuit needed to compute an answer. Here the input would be encoded in $n$ qubits and lives in the Hilbert space $\ell^2(\{0,1\}^n)$, and a simple gate in the circuit is a reversible unitary transformation on some of the qubits. The output is a measurement of some qubits. Quantum mechanics does not allow us to perform exact measurements, and so our output is naturally probabilistic. A quantum  algorithm solves a decision problem, iff the probability that it outputs the correct answer (Yes/No) is $\geq \frac23$ (this constant can be changed).   The quantum analogues of $\P$ and $\NP$ are $\BQP$ and $\QMA$: $\BQP$ is the class of decision problems for which there is a polynomially-sized quantum circuit computing the answer with high probability, and $\QMA$ is the class of problems, for which when the answer is Yes, there exists a poly-sized quantum circuit verifying the answer with high-probability. The counting analogue of $\SP$ is thus $\#\BQP$ and can be thought of as counting the number of accepting witnesses to a $\QMA$ verifier.

\subsection{Algebraic Complexity Theory}


Arithmetic (algebraic) Complexity theory is the study of computing polynomials $f(x_1,\ldots,x_n)\in \mathbb{F}[x_1,\ldots,x_n]$ in $n$ formal variables using simple operations $*,+,-,/$, where the input are the variables $x_1,\ldots,x_n$ and arbitrary constants from the underlying field. The complexity of the given polynomial $f$ is then the minimal number of such operations needed to compute the polynomial within the given model. There are three basic models of computations -- formulas, algebraic branching programs (ABPs) and circuits. For details on Algebraic Complexity theory see~\cite{BCS,Bur00}. Throughout the polynomials $f$ will be assumed to have $O(poly(n))$ bounded total degrees.

The algebraic complexity classes $\VP$ and $\VNP$ were introduced by Valiant~\cite{V1,V2}, as the algebraic analogues of $\P$ and $\NP$ (we refer to~\cite{Bur00} for formal definitions and properties).  

The \defn{class $\VP$} is the class of sequences of polynomials for which there is a constant $k$ and a $O(n^k)$-sized arithmetic circuit computing them. By \defn{arithmetic circuit} we mean a directed acyclic graph with source nodes containing  variables $x_1,\ldots,x_n$ or constants from the field, and the other vertices contain one simple operation performed with input from the nodes pointing to that vertex.  There is only one sink, which should contain the result of all the computations, our polynomial $f$. Let $S(f)$ denote the minimal possible size of a circuit computing $f$.

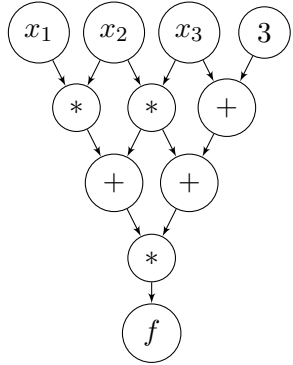
\begin{ex}

\begin{tabular}{p{2in}p{3in}}

\,

\begin{tikzpicture}

\tikzset{vertex/.style = {shape=circle,draw}}
\tikzset{edge/.style = {->,> = latex'}}
\node[vertex] (a) at  (0,4) {$x_1$};
\node[vertex] (b) at  (1,4) {$x_2$};
\node[vertex] (c) at  (2,4) {$x_3$};
\node[vertex] (d) at  (3,4) {$3$};
\node[vertex] (a1) at  (0.5,3) {$*$};
\node[vertex] (b1) at (1.5,3) {$*$};
\node[vertex] (c1) at (2.5,3) {$+$};
\node[vertex] (a2) at  (1,2) {$+$};
\node[vertex] (b2) at (2,2) {$+$};
\node[vertex] (a3) at (1.5,1) {$*$};
\node[vertex] (a4) at (1.5,0) {$f$};

\draw[edge] (a) to (a1);
\draw[edge] (b) to (a1);
\draw[edge] (b) to (b1);
\draw[edge] (c) to (b1);
\draw[edge] (c) to (c1);
\draw[edge] (d) to (c1);
\draw[edge] (a1) to (a2);
\draw[edge] (b1) to (a2);
\draw[edge] (b1) to (b2);
\draw[edge] (c1) to (b2);
\draw[edge] (a2) to (a3);
\draw[edge] (b2) to (a3);
\draw[edge] (a3) to (a4);

\end{tikzpicture} &
\

This circuit computes the polynomial $f= x_2x_3(x_1+x_2)(3+x_3)$ using 4 input nodes and 6 internal operations.
\end{tabular}
\end{ex}

\defn{The class $\VPs$} is the class of polynomials $f$, which have $O(n^k)$-sized formulas. A \defn{formula} is a circuit whose graph is a binary tree, so no output can be used twice. Let $L(f)$ denote the minimal formula size of $f$. Then a formula is recursively composed of operations  $f = g*h$ or $f=g+h$, so we have $L(f)\leq L(g) + L(h) + 1$.

\begin{ex}
Let $f = x_1^2 +x_2^2+x_3^2 + x_1x_3 +x_1x_2+3x_2x_3 = (x_1+x_2)*(x_1+x_3) + (x_2+x_3)(x_2+x_3)$, which has formula size $3+1+3$, so $L(f) \leq 7$.  
\end{ex}

We have that $S(f) \leq L(f)$ by definition, and according to~\cite{VSBR}, $S(f) \leq L(f)^{\log n}$. 

Finally, \defn{the class $\VBP$} is the class of polynomials $f$ which can be computed with a poly-sized Algebraic Branching Program. Informally, this is a directed acyclic graph with $\deg(f)$-many layers, one source $s$ and one sink $t$, each edge $e$ is labeled by a linear function in the variables $x_1,\ldots,x_n$ called $w(e)$ and the output is computed by going over all directed paths $p$: 
$$f = \sum_{p:s\to t} \, \prod_{e \in p} w(e).$$
The size of the branching program is defined as the maximal number of nodes in a layer, and by our assumptions is polynomially equivalent to the size of the given graph. Let $M(f)$ be the minimal size ABP needed to compute $f$, then $\VBP$ is the class of families of polynomials $f_n$ for which there is a fixed $k$ with $M(f_n)=O(n^k)$.

\defn{The class $\VNP$} is the class of polynomials $f$, such that there exists a fixed $k$ and a polynomial $g \in \VP$ with  $m=O( n^k )$ variables, such that
$$f(x_1,\ldots,x_n) = \sum_{b \in \{0,1\}^{m-n}} g(x_1,\ldots,x_n,b_1,\ldots,b_m).$$
In particular, every polynomial whose coefficients in the monomial expansion are easy to determine, would be in $\VNP$.

It is clear that $\VPs \subset\VBP \subset \VP \subset \VNP$, but are these classes different?

\begin{conj}[Valiant]
We have that $\VBP \neq \VNP$.
\end{conj}

We also believe that $\VP \neq \VNP$, but this problem is even harder to approach.


 As Valiant showed, for every polynomial $f \in \mathbb{C}[x_1,\ldots,x_n]$  there exists a $K$ and a $K \times K$ matrix $A$, s.t. $A = A_0 + \sum_{i=1}^N A_i x_i$ with $A_j \in \mathbb{C}^{K \times K}$,  such that 
 $$\det A= f.$$  
 The smallest such $K$ is the \defn{determinantal complexity $\dc(f)$} of $f$ and it is  \emph{finite for every $f$.} 
 
\begin{ex}
 Let $f = x_1^2 + x_1x_2 +x_2x_3 +2x_1$, then let
 $$A: =  \begin{bmatrix}
x_1 +2  & x_2 \\ -x_3+2 & x_1+x_2 
\end{bmatrix},$$
so $f=\det A$. Since $\deg(f)=2$, the smallest such matrix would be $2\times 2$ and so  $\dc(f) = 2$.
\end{ex} 

As Valiant also showed, we have
$$\dc(f) \leq 2 L(f),$$
and also that $M(f)$ and $\dc(f)$ are polynomially equivalent.
Thus from now on we can use  \emph{ $\dc$ as complexity measure}. We say that \defn{determinant is universal for $\VBP$}, in the sense that $f \in \VBP$ iff $\dc(f) = poly(n)$.

The classical universal $\VNP$-complete polynomial is the \defn{permanent}
$$\per_m [X_{ij}]_{i,j=1}^m = \sum_{\sigma \in S_m} \prod_{i=1}^m X_{i\sigma(i)}$$
in the sense that every $f \in \VNP$ of $\deg(f)\leq n^c$ can be written as a $\per_m[A]$ for some matrix $A$ of affine linear entries, and of polynomial size $m=O(n^k)$. It is much more powerful than the determinant.

Thus, to show that $\VBP \neq \VNP$ we need to show the following

\begin{conj}[Valiant~\cite{V2}]\label{conj:dc}
The determinantal complexity $\dc(\per_m)$ grows
superpolynomially in~$m$.
\end{conj}

It is known that $\dc(\per_m) \le 2^m -1$~\cite{Gre}, and $\dc(\per_m) \geq \frac{m^2}{2}$~\cite{MR}.

\medskip

The connection between $\P$ vs $\NP$ and $\VP$  vs $\VNP$ is exemplified in the following statement.

\begin{theorem}[\cite{Bur0}]
If one shows that $\VP=\VNP$ over a finite field then $\P=\NP$. If $\VP=\VNP$ over $\mathbb{C}$ and the Generalized Riemann Hypothesis holds then $\P=\NP$.
\end{theorem}

From here on we will only work over the field of constants $\mathbb{C}$.

We believe that separating $\VBP$  from $\VNP$ is the easier problem as the algebraic structure gives more tools. An approach to show $\VBP \neq \VNP$ is the \emph{Geometric Complexity Theory}, which will be discussed in Section~\ref{s:ac_in_gct}.

\section{Applications of Computational Complexity in Algebraic Combinatorics}\label{s:cc_in_ac}

Here we discuss some open problems in Algebraic Combinatorics. These can be phrased more formally using computational complexity theory and potentially answered within its framework.

\subsection{Open problems: combinatorial interpretation}

Semistandard Young Tableaux, the hook-length formula, the RSK correspondence, the Littlewood-Richardson rule are all examples of beautiful combinatorics. Besides aesthetically appealing, such results are also quite useful. Within Representation Theory they provide effective tools to understand the structure of group representations. Within asymptotics and statistical mechanics they give tools to understand behavior of lozenge tilings (dimer covers of the hexagonal grid), longest increasing subsequences of permutations, behavior of random sorting networks, random matrix eigenvalues etc.

Following the discovery of the Littlewood-Richardson rule in 1934, Murnaghan~\cite{Mur38} defined the Kronecker coefficients of $S_N$ and observed that computing even simple special cases is difficult. Interest in specifically nonnegative combinatorial interpretation can be found in~\cite{Las79,GR85}, and was formulated explicitly by Stanley as Problem 10 in his list ``Open Problems in Algebraic Combinatorics''~\cite{Sta00}\footnote{See this for the original list and updates on the problems \url{https://mathoverflow.net/questions/349406/}}.

\begin{op}\label{op:kron}
Find a combinatorial interpretation of $g(\la,\mu,\nu)$, i.e. a family of combinatorial objects $C(\la,\mu,\nu)$, such that $g(\la,\mu,\nu)=|C(\la,\mu,\nu)|$. 
\end{op}

Over the years, there has been very little progress on the question. In 1989 Remmel determined $g(\la,\mu,\nu)$ when two of the partitions are hooks~
\cite{R89}. In 1994 Remmel and Whitehead~\cite{RW} determined $g(\la,\mu,\nu)$ when two of the partitions are two-rows, i.e. $\ell(\la), \ell(\mu)\leq 2$. This case was subsequently studied also in~\cite{BMS}. In 2006 Ballantine and Orellana~\cite{BO} determined a rule for $g(\la,\mu,\nu)$ when one partition is a two-row, e.g. $\mu=(n-k,k)$, and the first row of one of the others is large, namely $\la_1 \geq 2k-1$. The most general rule was determined by Blasiak in 2012~\cite{Bla} when one partition is a hook, and this was later simplified by Blasiak and Liu~\cite{BL,Liu17}; informally it states that $g(\la,\mu,(n-k,1^k))$ is equal to the number of tableau in $\bar{1}<1<\bar{2}<\cdots$ of shape $\la$, type $\mu$ with restrictions on the ordering and certain entries. Other very special cases have been computed in the works of Bessenrodt-Bowman~\cite{BB17} for multiplicity-free products; when the marginals correspond to pyramids in Ikenmeyer-Mulmuley-Walter~\cite{IMW17}; near-rectangular partitions by Tewari~\cite{Tew15} etc. 

\begin{rem}
Most combinatorial interpretations in the area count tableaux or permutations with various restrictions. That, however, should not limit our scope. Consider the following labeled rooted \emph{partition trees}
 $T(m,\ell,k)$ whose vertices are labelled by $(a,b,\lambda,j)$, $j \leq ab$, $\lambda \vdash b$. The leaves correspond to labels with  $b=1$ and can thus be labeled by only  $(a, j)$ with $0\leq j \leq a$. 
Let the root be  $(m,\ell,\lambda,k)$ for some $\lambda \vdash \ell$. 

We impose the following local conditions between vertices and their children.
Let a vertex be labeled $(a,b,\lambda,j)$, with $\lambda = (1^{b_1},\ldots,n^{b_n})$. Then it has   at most $n$ children and their labels are of the form  $(a_1,b_1,\lambda^1,j_1),\ldots, (a_n,b_n,\lambda^n,j_n)$. s.t.
\begin{itemize}
\item $a_i = i(a+2)- 2(\lambda'_1+\cdots+\lambda'_i)$ for all $i=1, \ldots,n$.
\item $j_1+\cdots+j_n = j - 2\sum_i (i-1)\la_i$. 
\end{itemize}
Finally, let  the leaves be  $\{(a_0,i_0),\ldots,(a_t,i_t)\}$. Then we must have for each $r<t$: $i_r \geq 2(i_{r+1}+\cdots+i_t) - ( a_{r+1}+\cdots+a_t)$.

\begin{thm}[Pak-Panova, 2014, see~\cite{Pporto}]
The Kronecker coefficient $g((m^\ell), (m^\ell), (m\ell-k,k))$ is equal to the number of partition trees $T(m,\ell,k)$.
\end{thm}

The proof follows from two observations. The first is the fact  that $g((m^\ell), (m^\ell), (m\ell-k,k)) = p_k(m,\ell)-p_{k-1}(,m,\ell)$, where $p_k(m,\ell)=\#|\{ \al \vdash k, \al \subset (m^\ell)\}|$ is the number of partitions of $k$ fitting inside a $m\times \ell$ rectangle  (see e.g.~\cite{PP14,MPPe,Val97}). Alternatively, these are the coefficients in the expansion of the $q$-binomials
$$\sum_k p_k(m,\ell)q^k = \binom{m+\ell}{m}_q.$$
The second part is to unwind the recursive proof of the unimodality of those coefficients via  Zeilberger's KOH identity~\cite{Z}. The recursion then gives the tree $T$. 
\end{rem}

\medskip

Motivated by other developments, further questions on the Kronecker coefficients have appeared. Following their work in~\cite{HSTZ} on the square of the Steinberg character for finite groups of Lie type, Saxl conjectured that $g(\delta_k,\delta_k,\mu)>0$ for all $k$ and $\mu \vdash \binom{k+1}{2}$ and $\delta_k=(k,k-1,\ldots,1)$ is the staircase partition. This was initially studied in~\cite{PPV}, where its generalization was formulated as
\begin{conj}[Tensor square conjecture]
For every $n\geq 9$ there exists a symmetric partition $\la \vdash n$, such that $\mathbb{S}_\la \otimes \mathbb{S}_\la$ contains every irreducible $S_n$ module. In other words, $g(\la,\la,\mu)>0$ for all $\mu\vdash n$.
\end{conj}

The above conjecture raises the question on simply determining when $g(\la,\mu,\nu)>0$. Advances on the tensor square conjecture were initially made in \cite{PPV,I15,LS}, see~\cite{Pan22} for a list of many more recent works. It is a consequence of representation theory that for $n>2$ for every $\mu\vdash n$ there is a $\la$, such that $g(\la,\la,\mu)>0$ (see ~\cite[Ex. 7.82]{S1}), but even that has no combinatorial proof.

Positivity results were proved using a combination of three methods -- the semigroup property constructing recursively positive triples from building blocks, explicit heighest weight constructions using the techniques in \cite{Ful}, and an unusual comparison with characters, which was originally stated in by Bessenrodt and Behns~\cite{BB} for when  $g(\la,\la,\la)>0$, later generalized in~\cite{PPV}, and in its final form in~\cite{PPq}.

\begin{prop}[\cite{PPq}]\label{prop:char_kron}
Let $\la, \mu \vdash n$ and $\la=\la'$. Let $\hat{\la} = (2\la_1-1,2\la_2-3,3\la_3-5,\ldots)$ be the principal hooks of $\la$. Then 
$$g(\la,\la,\mu)\geq |\chi^\mu(\hat{\la})|.$$
\end{prop}

In 2020, with Christine Bessenrodt we generalized the conjecture as follows.

\begin{conj}[Bessenrodt-Panova 2020]
For every $n$ there exists a $k(n)$, such that for every $\la \vdash n$ with $\la=\la'$ and $d(\la)>k_n$ which is not the square partition, we have $g(\la,\la,\mu)>0$ for all $\mu \vdash n$.
\end{conj}

Here $d(\la) = \max\{i: \la_i\geq i\}$ is the Durfee square size of the partition. Partial progress on that conjecture will appear in the work of Chenchen Zhao.

\medskip

Another question motivated by Quantum Information Theory, pertaining to the so called ``\emph{quantum marginal problem}'', is for which triples of rational vectors $(\al,\be,\ga)$ is $g( k \al, k\be, k\ga)>0$ for some $k \in \mathbb{N}$, see e.g.~\cite{CHM}. Thanks to the semigroup property these triples form a convex polytope, a special case of the so-called \emph{Moment polytope} of a compact connected Lie group and a unitary representation.  The Kronecker polytope can actually be described in a certain concrete sense, see~\cite{VW}. 
The analogous question on positivity of Littlewood-Richardson coefficients is the Horn problem on spectra of matrices $A, B, C$, such that $A+B=C$. The resolution of the ``Saturation conjecture'' in~\cite{KT99} established that the inequalities cutting out the polytope of eigenvalue triples coincide with the inequalities defining triples of positive LR coefficients.

\medskip

Similar questions pertain to the plethysm coefficients. The following problem is number 9 in Stanley's list~\cite{Sta00}.

\begin{op}[Stanley]
Find a combinatorial interpretation for the plethysm coefficients $a^\la_{\mu,\nu}$.
\end{op}

Even the simple case for $a_\la(d[n]) = \< s_\la , h_d[h_n] \>$ is not known. 

A detailed survey on the partial results and methods can be found in~\cite{COSSZ}. 

There is no direct connection between the Kronecker and plethysm coefficients. Yet we know that when $\ell(\la)\leq 2$
$$a_\la(d[n]) = g(\la,n^d,n^d).$$
 An inequality between them in their stable limit is given in Theorem~\ref{thm:kron_pleth} and was obtained using their interpretations within GCT.

There is one major conjecture on plethysm coefficients.

\begin{conj}[Foulkes]\label{conj:foulkes}
Let $d > n$, then 
$$a_\la(d[n]) \geq a_\la(n[d])$$
for all $\la \vdash nd$. 
\end{conj}

This conjecture is related to Alon-Tarsi conjecture, and appeared in the GCT approaches as well. In~\cite{DIP} we proved it for some families of 3-row partitions, see Section~\ref{ss:mult_obstr}. 

\subsection{Complexity problems in Algebraic Combinatorics}\label{ss:comp_prob}

We will now study the important quantities in Algebraic Combinatorics with respect to their computational complexity leading to a classification by such ``hardness''. This gives a paradigm to understand these constants and either explain when a nice formula could be found, or  show that a combinatorial interpretation is unlikely as it would violate computational complexity assumptions like $\P \neq \NP$. 

Such questions have been quite common in other branches of combinatorics, like Graph Theory, and many of the graph theoretic problems are at the heart of Computational Complexity. Investigating computational complexity properties of structure constants was initiated when Algebraic Complexity Theory was developed. It came to prominence when Geometric Complexity Theory put understanding Littlewood-Richardson, Kronecker and plethysm coefficients in the spotlight. Most recently, understanding computational complexity  has been developed as a framework to formalize combinatorial properties of its own interest as in~\cite{Pak22}. 

\begin{ex}
Consider the problem \#SYT: given the input partition $\lambda$, compute the number of standard Young tableaux of shape $\la$. The answer would depend on how the input is encoded. Suppose that $\lambda$ is encoded in unary, i.e. each part $\la_i$ takes up $\la_i$ bits, and so the input size is $n=|\la|$. Using the HLF formula we can compute the product in $O(n)$ time and thus the problem is in $\FP$. If the input is in binary, then the input size is $|I| = O( \log_2(\la_1) \ell(\la))$ and $n=|\la| = O(2^{|I|})$. For most such partitions we would have $f^\la =o(2^n) =o(2^{2^{|I|}})$. This answer is too big, as it would require exponential space to even encode. This shows that binary input would not be appropriate for this problem at all. 
\end{ex}

As the example shows, the number of SYTs of shape $\la$ can be computed in polynomial time (when the input is in unary). We have that $f^\la = K_{\la,1^n}$, so the next natural problem is to compute the number of  SSYTs of shape $\la$ and  given type $\al$. This time, however, there is no product nor determinantal formula known.

\smallskip

\textsf{KostkaPos}: 

\texttt{Input:} $(\la_1,\la_2,\ldots),(\al_1,\al_2,\ldots)$

\texttt{Output:} Is $K_{\la,\al}>0$? 

\smallskip

This is the problem on deciding positivity of Kostka numbers. We know that $K_{\la,\mu}>0$ if and only if $\la \succ \mu$ in the dominance order, which is the set of linear inequalities for every $i=1,\ldots,\ell(\la)$
$$\la_1+\cdots +\la_i \geq \mu_1+\cdots +\mu_i.$$
Thus, given $\la$ and $\al$ in either binary or unary, we can check these inequalities in $O(\ell(\la))$ time,  so 
$  \textsc{KostkaPos} \in \P. $

The computational problem, however, is far from trivial

\smallskip

\textsf{ComputeKostka}: 

\texttt{Input:} $(\la_1,\la_2,\ldots),(\al_1,\al_2,\ldots)$ 

\texttt{Output:} Value of $K_{\la,\al}$.

\smallskip

\begin{thm}[Narayanan \cite{Nar}]
When the input $\la,\al$ is in binary, \textsf{ComputeKostka} is $\SP$-complete.
\end{thm}

It is not apriori clear why the problem (with binary input) would be in $\SP$ given that the SSYTs themselves have exponentially many entries. Yet $K_{\la,\al}$ can be computed as the number of integer points in the Gelfand-Tsetlin polytope, defined by $O(\ell(\al)^2)$ many linear inequalities. These inequalities can be verified in polynomial time. 
The proof of completeness uses a reduction to \textsf{KNAPSACK}, which is well known to be $\SP$-complete in binary, but it can be solved by a pseudopolynomial dynamic algorithm, so in unary it is in $\FP$. \footnote{The input encoding actually changes the problem and is usually part of the problem specification. In the early days of the development of Computational Complexity, a problem solvable in polynomial time when the input is in unary was said to be solvable in ``pseudopolynomial time''. When it is $\SP$-hard when the input is in unary, it would be called ``strongly $\SP$-hard'', see e.g.~\cite{GJ}.}
 
Yet, when the input is in unary for Kostka, in the general case, reduction to \textsf{KNAPSACK} does not give anything. Nonetheless, we conjecture that it is still hard. 

\begin{conj}[Pak-Panova 2020]
When the input is unary we have that \textsc{ComputeKostka} is \SP-complete.
\end{conj}

Here it is easy to see that the problem is in $\SP$, but not that it is hard. 

\medskip

Next we turn towards the Littlewood-Richardson coefficients. 

\smallskip

\textsf{LRPos}: 

\texttt{Input:} $\la,\mu, \nu$

\texttt{Output:} Is $c^{\la}_{\mu\nu}>0$? 

\smallskip

The proof of the Saturation Conjecture by Knutson and Tao~\cite{KT99} showed that an LR coefficient is nonzero if and only if the corresponding hive polytope is nonempty, see~\cite{DM06,MNS12}. This polytope is a refinement of the Gelfand-Tsetlin polytope, defined by $O(\ell(\la)^3)$ many inequalities. Showing that the polytope is nonempty is thus a linear programming problem, which can be solved in polynomial time. Thus
\begin{thm}
We have that \textsc{LRPos} $\in \P$ when the input is binary (and unary ditto). 
\end{thm}

\textsf{ComputeLR}: 

\texttt{Input:} $\la,\mu, \nu$ 

\texttt{Output:} Value of $c^{\la}_{\mu\nu}$.

Using the polytopal interpretation to show that even when the input is in binary we have \textsf{ComputeLR}$\in \SP$, the fact that Kostka is a special case of LR, see~\eqref{eq:kostka-lr}, we get the following.

\begin{thm}[Narayanan \cite{Nar}]
When the input $\la,\al$ is in binary, \textsf{ComputeLR} is \SP-complete.
\end{thm}

Yet again, when the input is in unary, we do not know whether the problem is still that hard.

\begin{conj}[Pak-Panova 2020]
When the input is in unary we have that \textsf{ComputeLR} is \SP-complete.
\end{conj}

We have that computing LR coefficients is in $\SP$  thanks to the Littlewood-Richardson rule and its polytopal equivalent formulation. If the input is unary, then the LR tableaux are the polynomial verifier, and one can check in $O(n^2)$ time if the tableaux satisfies all the conditions. The hard part here again is to show that computing them is still hard, namely that an $\SP$-complete problem like 3SAT would reduce to \textsf{ComputeLR}.

\medskip

None of the above has been possible for the Kronecker and plethysm coefficients, however, due to the lack of any positive combinatorial formula. 

\smallskip

\textsf{KronPos}: 

\texttt{Input:} $\la,\mu, \nu$

\texttt{Output:} Is $g(\la,\mu,\nu)>0$? 

\smallskip

The Kronecker coefficients have particular significance in GCT, see Section~\ref{s:ac_in_gct}. In the early stages Mulmuley conjectured~\cite{MS2} that they would be like the Littlewood-Richardson, so \textsf{KronPos}$\in \P$, which was recently disproved.

\begin{thm}[\cite{IMW17}]
When the input $\la,\mu,\nu$ is in unary, \textsf{KronPos} is $\NP$-hard.
\end{thm}

The proof uses the fact that in certain cases $g(\la,\mu,\nu)$ is equal to the number of pyramids with marginals $\la,\mu,\nu$, see~\cite{Val99}, and deciding if there is such a pyramid is $\NP$-complete. However, the problem is not yet in $\NP$, because we do not have  polynomially verifiable witnesses showing that $g(\la,\mu,\nu)>0$ when this happens. 

Needless to say, the problem would be even harder when the input is in binary, and we do not consider that here.

Mulmuley also conjectured that computing the Kronecker coefficients would be in $\SP$, again mimicking the Littlewood-Richardson coefficients.  
\smallskip

\textsf{ComputeKron}: 

\texttt{Input:} $\la,\mu, \nu$ 

\texttt{Output:} Value of $g(\la,\mu,\nu)$.

\smallskip

\begin{op}[Pak]
Show that \textsc{ComputeKron} is not in $\SP$ under reasonable complexity theoretic assumptions such as $\PH$ not collapsing. 
\end{op}

If the above is proven, that would make any solution to Open Problem~\ref{op:kron} as unlikely as the polynomial hierarchy collapsing. Any reasonable combinatorial interpretation as counting certain objects would show that the problem is in $\SP$, as the objects would likely be verifiable in polynomial time.  

Note that \textsf{ComputeKron}$\, \in\GapP$ (\cite{BI08}) as it is easy to write an alternating sum for its computation, for example using contingency arrays, see~\cite{PP17}. This further shows that $\SP$ would be a natural class for this problem as it is already in $\GapP_{\geq 0}$. 

The author's experience with Kronecker coefficients seems to suggest that some particular families would be as hard as the general problem. 

\begin{conj}[Panova]
We have that \textsf{ComputeKron} is in $\SP$  when $\ell(\la)=2$ if and only if \textsf{ComputeKron} is in $\SP$ in the general case. Likewise, \textsf{ComputeKron}   for $\mu=\nu =(n^d)$ and $\la \vdash nd$ as the input is in $\SP$ if and only if the general problem is in $\SP$.
\end{conj}

The last part concerns the \emph{rectangular Kronecker coefficients} of special significance in GCT, see Section~\ref{s:ac_in_gct} and~\cite{IP}. 

It is worth noting that when the partitions have fixed lengths, we have that \textsc{ComputeKron} is in $\FP$ even when the input is in binary, see~\cite{CDW,PP17}. Moreover, from the asymptotic analysis and symmetric function identities  in~\cite{PPd}, it follows that
\begin{prop}
Let $k$ be fixed and $(\la,\mu,\nu) \vdash n$ be partitions with diagonals at most $k$, i.e. $d(\la), d(\mu),d(\nu)\leq k$. Then $g(\la,\mu,\nu)$ can be computed in time $O(n^{4k^3})$. 
\end{prop}
Note that in this case we have $g(\la,\mu,\nu) \leq C_k n^{4k^3 + 13k^2 +31 k}$ for an explicit constant $C_k$. This in itself does not guarantee the efficiency of the algorithm computing them, but in this case can be easily derived. On the other hand,  when the lengths of the partitions are bounded by $k$  the efficient algorithms run in time $O( (\log n)^{k^3 \log k} )$. We do not expect that a similar efficient algorithm exists in the more general  case of fixed diagonals.

\bigskip

\begin{tabular}{p{2.5in}p{3in}}
\textsf{PlethPos}: 

\texttt{Input:} $\la,\mu, \nu$

\texttt{Output:} Is $a^\la_{\mu,\nu}>0$? 
&

\textsf{ComputePleth}: 

\texttt{Input:} $\la,\mu, \nu$

\texttt{Output:} Value of $a^\la_{\mu,\nu}$.

\end{tabular}

Using symmetric function identities, it is not hard to find an alternating formula for the plethysms and show that they are also in $\GapP$, see~\cite{FI20}. They also show that \textsf{PlethPos} is $\NP$-hard. We suspect that \textsf{ComputePleth} may not be in $\SP$ in the general case, but also when $\mu,\nu$ are single row partitions. The coefficient then $a_\la(d[n])$ has special significance in GCT, see Section~\ref{s:ac_in_gct}. 

\begin{op}
Determine whether \textsc{PlethPos}$\in \NP$ and \textsf{ComputePleth}$\in \SP$ under reasonable complexity theoretic assumptions.
\end{op}

\medskip

The representation theoretic significance of these structure constants poses the natural question on their computation via quantum circuits. Quantum computing can be powerful on algebraic and number theoretic problems. The structure constants in question are dimensions of various vector spaces, and it is natural to expect that such quantities could be amenable to efficient quantum computation. 
While this is not straightforward, Beals' quantum Fourier transform over the symmetric group~\cite{B97} gives the following
\begin{thm}
\textsc{KronPos} is in $\QMA$. \textsc{ComputeKron} is in $\#\BQP$.
\end{thm}

These statements have been claimed in~\cite{HCW}.  The first statement and a weaker version of the second were shown in~\cite{BCGHZ}. The full proof of the second statement appears in~\cite{IS}.
As the first group noted, the statements should be true even when the input is in binary. 

\begin{op}
Show that when the input $(\la,\mu,\nu)$ is in binary, then \textsc{KronPos} is in $\QMA$ and \textsc{ComputeKron} is in $\BQP$.
\end{op}
 
With the input in binary we can no longer use the symmetric group $S_n$ as $n$ would be too large and we will have to use the $GL$ interpretation of the Kronecker coefficients.

\subsection{Proof of concept: character squares are not in \SP}

Underlying all the representation theoretic multiplicities mentioned above are the characters of the symmetric group. For example, equation~\eqref{eq:char_kron} expresses the Kronecker coefficients via characters, and the other structure constants can also be expressed in similar ways. What then can we say about computing the characters and can this be used in any way to help with the problems in Section~\ref{ss:comp_prob}

The characters satisfy some particularly nice identities coming from the orthogonality of the rows and columns of the character table in $S_n$.  We have that 
\begin{equation}\label{eq:char_sum}
\sum_{\la \vdash n} \chi^\la(w)^2 = \prod_i i^{c_i} c_i! \, ,
\end{equation}
where $c_i =$ number of cycles of length $i$ in $w \in S_n$. When $w=\id$, we have that $\chi^\la(\id) = f^\la$, the number of SYTs and the identity becomes equation~\ref{eq:rsk}. That equation, as mentioned in Section~\ref{ss:bakcground}, can be proven via the beautiful RSK bijection. The first step in this proof is to identify $(f^\la)^2$ as the number of pairs of SYTs of the same shape.

Could anything like that be done for equation~\eqref{eq:char_sum}? The first step would be to understand what objects $\chi^\la(w)^2$ counts, does it have any positive combinatorial interpretation? We formulate it again using the CC paradigm as

\smallskip
\textsf{ComputeCharSq}: 

\texttt{Input}: $\la,\al \vdash n$, unary. 

\texttt{Output}: the integer $\chi^\la(\al)^2$.

\smallskip

\begin{thm}[\cite{IPP22}]
\textsf{ComputeCharSq}$\not \in \SP$ unless $PH=\Sigma^{\P}_2$. 
\end{thm}

The last condition says ``polynomial hierarchy collapses to the second level'', which is almost as unlikely as $\P=\NP$, and is a widely believed complexity theoretic assumption. The proof uses the intermediate vanishing decision problem

\smallskip
\textsf{CharVanish}: 

\texttt{Input}: $\la,\al \vdash n$, unary. 

\texttt{Output}: Is  $\chi^\la(\al)=0?$

\smallskip
\begin{thm}[\cite{IPP22}]\label{thm:char_vanish}
We have that \textsf{CharVanish} is \CeqP-complete under many-to-one reductions.
\end{thm}

In order to prove this we use the Jacobi-Trudi identity to write $\chi^\la(\al)$ as an alternating sum of ordered set partition functions. 
Let $\la \vdash n$ with $\ell(\la)\leq \ell$, and let $\al$ be a composition of $n$. Then
$$\chi^\la(\al) \. = \. \sum_{\sigma \in S_{\ell(\la)}} \sgn(\sigma) \. P(\al,\la+\sigma-\id).$$
Using number theoretic restrictions we limit the entries to just two:

\begin{prop}\label{p:char-part3}
Let $\mathbf{c}$ and $\mathbf{d}$ be two sequences of nonnegative integers, such that $|\mathbf{c}|=|\mathbf{d}|+6$. Then there are partitions $\la$ and $\al$ of size $O(\ell |\mathbf{c}|)$ determined in linear time, such that
$$\chi^\la(\al) \. = \. P\big(\mathbf{c}, \ov{\mathbf{d}}\big) \. - \. P\big(\mathbf{c},\ov{\mathbf{d}'}\big),$$
where  $\ov{\mathbf{d}}:=(2,4,d_1,d_2,\ldots)$ and $\ov{\mathbf{d}'} := (1,5,d_1,d_2,\ldots)$.
\end{prop}

We then use the fact that matchings can be encoded as set partition problems, by encoding the edges/hyperedges as unique integers in a large basis as in~\cite{GJ}. After some constructions, putting two 3d-matching problem instances on one hyepergraph, with hyperedges of 4 vertices we conclude that

\begin{prop}[\cite{IPP22}]
For every two independent 3d matching problem instances $E$ and $E'$, there exist $\mathbf{c}$ and $\mathbf{d}$ as above, such that 
$$\#3DM(E)- \#3DM(E') = \frac{1}{\delta}\left(  P\big(\mathbf{c}, \ov{\mathbf{d}}\big) \. - \. P\big(\mathbf{c},\ov{\mathbf{d}'}\big) \right) = \frac{1}{\delta} \chi^\la(\al),$$
where $\delta$ is a fixed multiplicity factor equal to the number of orderings. 
\end{prop}

Finally, we observe that counting 3d matchings is a $\SP$-complete problem. Thus the last equations shows that $\chi^\la(\al)=0$ iff $\#3DM(E)= \#3DM(E')$, i.e. vanishing is equivalent to two independent $\SP$ functions being equal. This makes it $\CeqP$-complete and proves Theorem~\ref{thm:char_vanish}.

To show the next steps we use classical CC results.
If $\chi^2 \in \SP$ then $[\chi^2 >0] \in \NP$, so $[\chi \neq 0] \in \NP$
and hence, by definition, $[\chi=0]\in \coNP$. Thus  $\CeqP \subset \coNP $. 
By a result of Tarui we have $\PH \subset \NP^{\CeqP}$, and from the above we  get $\PH \subset \NP^{\coNP}=\Sigma^\P_2$. So $\PH = \Sigma^\P_2$, and the proof follows. 

\medskip

In contrast with this result, we note that Beals' quantum Fourier transform over the symmetric group~\cite{B97} actually gives an efficient quantum algorithm for the characters.

\section{Applications of Algebraic Combinatorics in Computational Complexity Theory}\label{s:ac_in_gct}

\subsection{Geometric Complexity Theory}\label{ss:gct}

Towards answering Conjecture~\ref{conj:dc} and showing that $\VBP \neq \VNP$, 
Mulmuley and Sohoni~\cite{MS1,MS2} proposed an approach based on algebraic geometry and
representation theory, for which they coined the name Geometric Complexity Theory (GCT). For an accessible and detailed treatment we refer to~\cite{BI}. 

Informally, the idea is to show that an $m \times m$ permanent of a variable matrix $[X_{i,j}]_{i,j=1}^m$  cannot be expressed as an $n\times n$ determinant of a matrix  with affine linear forms as entries for $n = O(m^k)$ for any $k$. Set $\mathbf{X}=(Z,X_{11},X_{12},\ldots,X_{mm})$ as the vector of variables in the matrix $X$ plus the variable $Z$ for the affine terms. Because we are considering all possible linear forms, we are looking at ${\det}_n[M \mathbf{X}]$ for all matrices $M \in \mathbb{C}^{n^2 \times n^2}$ and we want to explore when $\per_m[X] = {\det}_n[M \mathbf{X}]$. Replacing these matrices by invertible ones, and then taking the Euclidean closure would give us a, slightly larger, space of polynomials containing $\{ {\det}_n[M \mathbf{X}] : M \in \mathbb{C}^{n^2 \times n^2} \} \subset \overline{ {\det}_n(GL_{n^2} \mathbf{X}) }$. Here the tools of Algebraic Geometry and Representation theory become available,  we can compare the orbit closures of $\per_m$ and ${\det}_n$ via their irreducible representations to show that containment is not possible for $n=poly(m)$.  

More formally, as outlined in~\cite{BIP}, the setup is as follows. Denote by $\Sym^n V^*$ the space of homogeneous polynomial functions of degree~$n$ on 
a finite dimensional complex vector space $V$.
The group $G:=\GL(V)$ acts on $\Sym^n V^*$ in the canonical way by linear substitution:
$(h \cdot f)(v) := f(h^{-1}v)$ for $h\in G$, $f\in \Sym^n V^*$, $v\in V$.  
We denote by $G\cdot f :=\{ hf \mid h\in G\}$ the \defn{ orbit} of~$f$.
We assume now $V:=\C^{n\times n}$, view the determinant $\det_n$ as an element of $\Sym^n V^*$,
and consider its \defn{ orbit closure}: 

\begin{equation}\label{def:Omegan}
\Det_n := \ol{\GL_{n^2}\cdot \det_n} \subseteq \Sym^n (\C^{n\times n})^* 
\end{equation}
with respect to the Euclidean topology, which is also the same as with respect to the Zariski topology.

For $n > m$ we consider the \defn{ padded permanent} defined as 
$X_{11}^{n-m} \per_m \in\Sym^n (\C^{m\times m})^*$ (here we replace the extra variable $Z$ mentioned in the beginning by $X_{11}$ directly). 
Via the standard projection $\C^{n\times n} \to \C^{m\times m}$, we can view 
$X_{11}^{n-m} \per_m$ as an element of the bigger space $\in\Sym^n (\C^{n\times n})^*$.

The following conjecture was stated in \cite{MS1}.

\begin{conjecture}[Mulmuley and Sohoni 2001]\label{conj:dc-bord}
For all $c\in\N_{\ge 1}$ 
we have $X_{11}^{m^c-m} \per_m \not\in \Det_{m^c}$ 
for infinitely many~$m$.
\end{conjecture}

As discussed in the beginning, if $\per_m = {\det_n}[M\mathbf{X}]$ for some $n$, using the fact that $\GL_{n^2}$ is dense in $\C^{n^2\times n^2}$, we have
that  $\dc(\per_m) \leq n$,
and $\per_m \in \Det_{n}$. 

Thus, Conjecture~\ref{conj:dc-bord} implies Conjecture~\ref{conj:dc}.

The following strategy towards Conjecture~\ref{conj:dc-bord} was proposed by Mulmuley and Sohoni in \cite{MS1}. 
We consider the space $\Det_n$ as an algebraic variety and study its structure via its coordinate ring.
Specifically, the action of the group $G=\GL(V)$ on 
$\Sym^n V^*$ induces an action on its graded coordinate ring 
$\C[\Sym^n V^*]=\oplus_{d\in\N} \Sym^d \Sym^n V$. 
The space $\Sym^d \Sym^n V$ decomposes into irreducible $GL_{n^2}$--modules with multiplicities exactly the \emph{plethysm} coefficients. 
The \defn{coordinate ring $\C[\Det_n]$ of the orbit closure $\Det_n$} 
is obtained as the homomorphic image
of $\C[\Sym^n V^*]$ via the restriction of regular functions, 
and the $G$-action descends on this. 
In particular, we obtain the degree~$d$ part $\C[\Det_n]_d$ of $\C[\Det_n]$
as the homomorphic $G$-equivariant image of $\Sym^d \Sym^n V$.

As a $G$--module, the coordinate ring $\C[\Det_n]$ is a direct sum of its 
irreducible submodules since $G$ is reductive. 
We say that $\la$ occurs in $\C[\Det_n]$ if 
it contains 
an irreducible $G$-module of type $\la$ and denote its multiplicity by $\delta_{\la,d,n}$, so we can write
\begin{equation}\label{eq:det_reps}
\IC[\Det_n]_d=\IC[\overline{GL_{n^2} {\det}_n}]_d \simeq \bigoplus_{\la \vdash nd} V_\la^{\oplus \delta_{\lambda,d,n}}
\end{equation}

On the other side, we repeat the construction for the permanent. Let \defn{$\ocp_{n,m}$ denote the orbit closure of the padded permanent} ($n>m$): 
\begin{equation}\label{def:Znm}
 \ocp_{n,m} := \ol{\GL_{n^2}\cdot X_{11}^{n-m} \per_m} \subseteq \Sym^n (\C^{n\times n})^*.
\end{equation}
If $X_{11}^{n-m} \per_m={\det}_n[M\mathbf{X}]$, then it is contained in $\Det_n$, then 
\begin{equation}\label{eq:per_in_det}
\ocp_{n,m} \subseteq \Det_n,
\end{equation} and the restriction defines 
a surjective $G$-equivariant homomorphism 
$\C[\Omega_n]\to\C[\ocp_{n,m}]$ of the coordinate rings. 
We can decompose this ring into irreducibles likewise, 
\begin{equation}\label{eq:per_reps}
\IC[\ocp_{n,m}]_d=\IC[\overline{GL_{n^2} \per_m^n}]_d \simeq \bigoplus_{\la \vdash nd} V_\la^{\oplus \gamma_{\lambda,d,n,m}}.
\end{equation}
If $\C[\Omega_n]\to\C[\ocp_{n,m}]$, then we must have $\ga_{\la,d,n,m}\leq \de_{\la,d,n}$ by Schur's lemma. 
A partition $\la$ for which the opposite holds, i.e.
\begin{equation}\label{eq:mult_obs}
\ga_{\la,d,n,m} > \de_{\la,d,n}
\end{equation}
is called a \defn{multiplicity obstruction}. Its existence shows that the containment~\eqref{eq:per_in_det} is not possible and hence the permanent is not an $n\times n$ determinant of affine linear forms. 
\begin{lemma}
If there exists an integer $d$ and a partition $\la \vdash n$, for which~\eqref{eq:mult_obs} holds, then $\dc(\per_m)>n$. 
\end{lemma}

The main conjecture in GCT is thus
\begin{conj}[GCT, Mulmuley and Sohoni~\cite{MS1}]\label{conj:mult-obstr}
There exist multiplicity obstructions showing that $\dc(\per_m) > m^c$ for every constant $c$. Namely, for every $n=O(m^c)$ there exists an integer $d$ and a partition $\la \vdash dn$, such that $\ga_{\la,d,n,m} > \de_{\la,d,n}$.
\end{conj}

A partition $\la$ which does not occur in $\IC[\Det_n]$, but occurs in $\IC[\ocp_{n,m}]$, i.e. $\ga_{\la}>0, \de_{\la}=0$, is called an 
{\em occurrence obstruction}. Its existence thus also
proves that $\ocp_{n,m} \not\subseteq \Det_n$ 
and hence $\dc(\per_m)>n$.

In \cite{MS1,MS2} it was suggested to prove Conjecture~\ref{conj:dc-bord} 
by exhibiting occurrence obstructions. 
More specifically, the following conjecture was stated. 

\begin{conjecture}[Mulmuley and Sohoni 2001]\label{conj:occ-obstr}
For all $c\in\N_{\ge 1}$, 
for infinitely many $m$, 
there exists a partition~$\la$ 
occurring in $\C[Z_{m^c,m}]$ but not in $\C[\Det_{m^c}]$.
\end{conjecture}

This conjecture implies Conjecture~\ref{conj:dc-bord}
by the above reasoning. 

\subsection{Structure constants in GCT}\label{ss:gct_str_const}

Conjecture~\ref{conj:mult-obstr} and the easier Conjecture~\ref{conj:occ-obstr} on the existence of occurrence obstructions  
has stimulated a lot of research and 
has been the main focus of researchers 
in geometric complexity theory.

Unfortunately, the easier Conjecture~\ref{conj:occ-obstr} turned out to be false.

\begin{theorem}[B\"urgisser-Ikenmeyer-Panova~\cite{BIP}]\label{thm:nooccobstr}
Let $n,d,m$ be positive integers with 
$n \ge m^{25}$ and $\la\vdash nd$. If $\la$ occurs in 
$\C[Z_{n,m}]$, then $\la$~also occurs in $\C[\Omega_n]$. 
In particular, Conjecture~\ref{conj:occ-obstr} is false. 
\end{theorem}

Before we explain its proof, we will establish the connection with Algebraic Combinatorics.

In \cite{MS1} it was realized that the GCT-coefficients $\ga_{\la,d,n}$
can be bounded by rectangular Kronecker coefficients,
we have $\ga_{\la,d,n}(\la) \le g(\la,n^d,n^d)$ for $\la\vdash nd$. 
In fact, the multiplicity of $\la$ in the coordinate ring of the orbit 
$\GL_{n^2}\cdot\det_n$ equals the so-called symmetric 
rectangular Kronecker coefficient ${\rm sk}(\la,n^d)$, see~\cite{BLMW}, which is in general defined as   
$${\rm sk}(\la,\mu):= \mult_\la \Sym^2(\mathbb{S}_\mu) \leq g(\la,\mu,\mu).$$ 

Note that an occurrence obstruction for 
$Z_{n,m}\not\subseteq\Omega_n$ 
could then be a partition $\la$ for which  $g(\la,n^d,n^d)=0$ 
and such that $\la$ occurs in $\C[Z_{n,m}]$. 
Since hardly anything was known about the actual coefficients $\ga_{\la,d,n}$, 
it was proposed in \cite{MS1} to find $\la$ 
for which the Kronecker coefficient $g(\la,n^d,n^d)$ vanishes 
and such that $\la$ occurs in $\C[Z_{n,m}]$.

\begin{conj}[\cite{MS1}]
There exist $\la$, s.t. $g(\la,n^d,n^d) =0$ and $\gamma_{\la,d,n,m}>0$ for some $n>poly(m)$. 
\end{conj}

This was the first conjecture to be disproved.

\begin{theorem}[\cite{IP}]\label{thm:no_kron}
Let {$n>3m^4$}, $\la \vdash nd$. If {$g(\la,n^d, n^d)=0$}, then $\gamma_{\la,d,n,m}=0$.
\end{theorem}

In order to show this, we need a characterization of $\gamma_{\la,d,n,m}$, which follows from the work of Kadish and Landsberg~\cite{KL}. 

\begin{prop}[\cite{KL}]
If $\gamma_{\la,d,n,m} >0$ then $\ell(\la)\leq m^2$ and $\la_1 \geq d(n-m)$.
\end{prop}

The rest revolves around showing that for such $\la$ the relevant Kronecker coefficients would actually be positive.

\begin{theorem}[\cite{IP}] \label{thm:kron_pos}
If 
$\ell(\la)\leq m^2$,
$\la_1 \geq nd - md$,
$d > 3 m^3$, and
$n > 3 m^4$, then $g(\la,n\times d, n \times d)>0$, except for 6 special cases.
\end{theorem}

The proof of this Theorem uses two basic tools. 

One is the result of Bessendrodt-Behns~\cite{BB}, generalized to Proposition~\ref{prop:char_kron}, that
$$g( k^k, k^k, k^k )>0$$
for all $k\geq 1$. 

The other is the semigroup property Proposition~\ref{thm:semigroup} applied in various combinations and settings together with conjugation. In particular, we also have that if $\al +_V \be:= (\al' + \be' )'$, the vertical addition of Young diagrams, we have $g(\al^1 + \be^1, \al^2 +_V \be^2, \al^3 +_V \be^3) >0$ whenever both $g(\al^1,\al^2,\al^3)>0$ and $g(\be^1,\be^2,\be^3)>0$.  

To prove the general positivity result, we cut our partitions $\la$ into squares of sizes $2\times 2, \ldots, m^2 \times m^2$ and some remaining partition $\rho$ with at most $m^3$ many columns bigger than 1, namely
$$\la = \sum_{k=2}^{m^2} n_k (k^k) + \rho.$$ The two rectangles can also be cut into such square pieces giving triples $(k^k, k^k,k^k)$ of positive Kronecker coefficients that can be combined together using the semigroup properties. Finally, for the remaining partition $\rho$, we show inductively that if $g(\mu, k_1^{k_1},k_1^{k_1})>0$ for some $\mu \vdash k_1^2$, then for all $k,\ell,r$ with $|\mu|+\ell+r=k^2$ we have all positive Kronecker coefficients
$$g\left( (\mu+\ell)+_V1^r, k^k, k^k\right) >0.$$

\medskip

Using the fact that determinantal complexity for all polynomials of fixed degree and number of variables  is finitely bounded, the GCT setup and the bounds on the multiplicities we obtain the following unexpected relation between rectangular Kronecker coefficients and plethysms. Note that the range of $d$ and $n$ here put these multiplicities in \emph{stable regime}, i.e. their values stabilize when $n,d$ increase.

\begin{theorem}[\cite{IP}]\label{thm:kron_pleth}
For every partition $\rho$, let $n \geq |\rho|$, $d\geq 2$, $\la := (nd-|\rho|,\rho)$. Then 
$$
g( \la , n^d, n^d) \geq a_\la(d[n]).
$$
\end{theorem}

 In fact, the proof gives $\sk(\la,n^d) \geq a_\la(d[n])$. 

\smallskip

The ideas in the proof of Theorem~\ref{thm:nooccobstr} are similar in philosophy, but technically different.  
We have that $\dc(X_1^s + \cdots +X_k^s) \leq ks$, as seen from the formula size relation and Valiant's proof~\cite{V1}. Then, after homogenization, we have $z^{n-s}(v_1^s+\cdots+v_k^s) \in \Omega_n$ for $n \geq ks$ and linear forms $z,v_1,\ldots,v_k$. 

Now we can consider 
\[
\Pow_{k}^n := \overline{\{\ell_1^n +  \cdots + \ell_k^n \mid \ell_i \in V \}} \in \Omega_{kn} ,
\]
and essentially prove, see also~\cite{DIP}, that using the same setup for coordinate rings replacing the determinant with the power sum polynomial, see Proposition~\ref{prop:dip_pow}, that
$${\rm mult}_\la(\IC[P_{k}^n]_d) = a_\la(d[n]) \text{ for }k\geq d$$
(for the partitions $\la$ of relevance).

 Comparing multiplicities then we get $\delta_{\la,d,n}=\mult_\la \mathbb{C}[\Omega_n] \geq a_{\la}(d[n])$. We show using explicit tableaux constructions, see~\cite{Ful}, that  $a_{\la}(d[n])>0$ for the 
 partitions $\la$ such that $\la_1 \geq d(n-m)$ and $\ell(\la) \leq m^2$. 
 
 \begin{rem}
 In~\cite{BIP} we show that occurrence obstructions don't work not just for permanent versus determinant, but for permanent versus power sum polynomial. Power sums are clearly much weaker computationally than the determinant polynomial. The barrier towards occurrence obstructions comes from the padding of the permanent, which results in partitions $\la$ with long first rows. The long first row makes the relevant multiplicities positive, as can be seen with the various applications of semigroup properties. 
 \end{rem}

\subsection{Multiplicity obstructions}\label{ss:mult_obstr}

In order to separate \VP from \VNP via determinant versus permanent it is now clear that occurrence obstructions would not be enough. To remedy this there are two approaches. 

We can replace the ${\det}_n$ by the Iterated Matrix Multiplication tensor ${\rm tr}(A_1 \cdots A_m)$, the trace of the product of $m$ matrices with affine linear entries of size $n\times n$. This is another $\VBP$ universal model, and the measure of complexity is $n$, the size of the matrices. In this case we will not be padding the permanent, and the partitions involved would not have long first rows. The drawback now is that computing the multiplicities is even more complicated.

Alternatively, we can look for \emph{multiplicity obstructions}, i.e. partitions $\la \vdash dn$, for which
$$ \gamma_{\la, d,n,m} < \delta_{\la,d,m} \text{ for some }n\gg poly(m),$$
where by $poly(m)$ we mean any fixed degree polynomial in $m$.

As a proof of concept, we consider another separation of polynomials, as done in~\cite{DIP}.

Consider the space $\IA_m^n := \IC[x_1,\ldots,x_m]_n$ of complex homogeneous polynomials of degree $n$ in $m$ variables.
Let $V := \IA_m^1$ be the space of homogeneous degree 1 polynomials.
We compare two subvarieties of $\IA_m^n$.
The first is the so-called \defn{Chow variety}
\[
\Ch_m^n := \{\ell_1 \cdots \ell_n \mid \ell_i \in V \} \subseteq \IA_m^n,
\]
which is the set of polynomials that can be written as a product of homogeneous linear forms.
In algebraic complexity theory this set is known as the set of polynomials that have homogeneous depth-two algebraic circuits of the form $\Pi^n \Sigma$, i.e., circuits that consists of an $n$-ary top product gate of linear combinations of variables.

The second variety is called a \defn{higher secant variety of the Veronese variety} and can be written as
\[
\Pow_{m,k}^n := \overline{\{\ell_1^n +  \cdots + \ell_k^n \mid \ell_i \in V \}} \subseteq \IA_m^n,
\]
which is the closure of the set of all sums of $k$ powers of homogeneous linear forms in $m$ variables, which also showed up in~\cite{BIP} as mentioned in \S\ref{ss:gct_str_const}.
In algebraic complexity theory this set is known as the set of polynomials that can be approximated arbitrarily closely by homogeneous depth-three powering circuits of the form $\Sigma^k \Lambda^n \Sigma$,
i.e., a $k$-ary sum of $n$-th powers of linear combinations of variables.

We now consider when $\Pow_{m,k}^n \not\subseteq \Ch_m^n$, or in other words, when is a power sum not factorizable  as a product of linear forms. While this is easy to see explicitly, here we will show
how GCT can work in practice when there are no \emph{occurrence obstructions}, namely, we will find \emph{multiplicity obstructions}.

The approach is in complete analogy to the approach described in Section~\ref{ss:gct} to separate group varieties arising from algebraic complexity theory. Here we've replaced $\per$ by a power sum polynomial, and $\det$ by the product of linear forms.

If $\Pow_{m,k}^n \subseteq \Ch_m^n$, then the restriction of functions gives a canonical $\GL_m$-equivariant surjection
\[
\IC[\Ch_m^n]_d \twoheadrightarrow \IC[\Pow_{m,k}^n]_d.
\]
Decomposing the two modules into irreducibles and comparing multiplicities for each $V_\la$ we have that
\begin{equation}\label{eq:obsineq}
\mult_\la(\IC[\Ch_m^n]_d) \geq \mult_\la(\IC[\Pow_{m,k}^n]_d).
\end{equation}
for all partitions $\la$ with $\ell(\la) \leq m$.
Therefore, a partition $\la$ that violates \eqref{eq:obsineq} proves that $\Pow_{m,k}^n \not\subseteq \Ch_m^n$ and is called a \emph{multiplicity obstruction}.
If additionally $\mult_\la(\IC[\Ch_m^n]_d)=0$, then $\la$ would be an \emph{occurrence obstruction}.

Since these are $GL_m$ modules we must have $\ell(\la)\leq m$   and since the total degree is $dn$ we have $\la \vdash dn$.

\begin{theorem}[\cite{DIP}]\label{thm:dip}
~\\\textup{(1) Asymptotic result: } Let $m\geq 3$, $n \geq 2$, $k=d=n+1$, $\la=(n^2-2,n,2)$. We have
$\mult_\la(\IC[\Ch_m^n]_d) < \mult_\la(\IC[\Pow_{m,k}^n]_d)$,
i.e., $\la$ is a multiplicity obstruction that shows $\Pow_{m,k}^n \not\subseteq \Ch_m^n$.

\noindent\textup{(2) Finite result: } In two finite settings we can show a slightly stronger separation:\\
\noindent \mbox{~}\ \ \textup{(a)} Let $k=4$, $n=6$, $m=3$, $d=7$, $\la=(n^2-2,n,2)=(34,6,2)$. Then
$\mult_\la(\IC[\Ch_m^n]_d) = 7 < 8 = \mult_\la(\IC[\Pow_{m,k}^n]_d)$, i.e., $\la$ is a multiplicity obstruction that shows $\Pow_{m,k}^n \not\subseteq \Ch_m^n$.\\
\noindent \mbox{~}\ \ \textup{(b)} Similarly, for $k=4$, $n = 7$, $m = 4$, $d = 8$, $\la=(n^2-2,n,2)=(47,7,2)$ we have $\mult_\la(\IC[\Ch_m^n]_d) < 11 = \mult_\la(\IC[\Pow_{m,k}^n]_d)$, i.e., $\la$ is a multiplicity obstruction that shows $\Pow_{m,k}^n \not\subseteq \Ch_m^n$.\\
\noindent Both separations \textup{(a)} and \textup{(b)} cannot be achieved using occurrence obstructions, even for arbitrary $k$:
for all partitions $\la$ of $\ell(\la)\leq m$ 
that satisfy $a_\la( d [n]) >0$
we have $\mult_\la(\IC[\Ch_m^n]_{d'})>0$ in these settings.
\end{theorem}

The proof involves two facts which relate the desired multiplicities with plethysms. We have that
$a_{\la}(d[n]) = \mult_\la(\IC[\A_m^n]_d)$ 
\begin{prop}[\cite{BIP}]\label{prop:dip_pow}
Let $\la\vdash dn$ with $\ell(\la) \leq m$. If $k\geq d$ then $\mult_\la \IC[\Pow_{m,k}^n]_d = a_\la(d[n])$.
\end{prop}
We also have that
\begin{lemma}[\cite{DIP}]
Let $\la \vdash nm$ with $\ell(\la)\leq m \leq n$. Then
$\mult_{\la}( \IC [\Ch_m^n]_d) \leq a_\la(n[d]).$
\end{lemma}
Finally we find explicit values and relations for the plethysm coefficients and prove in particular the following
\begin{thm}[\cite{DIP}]
Let $\la = (n^2-2,n,2) \vdash n(n+1)$ and let $d=n+1$. Then
$$a_\la(d[n]) =1 + a_\la(n[d]).$$ 
\end{thm}
In particular, this confirms Foulkes conjecture~\ref{conj:foulkes}.

\section{Discussion}

As we have seen, structure constants from Algebraic Combinatorics, mostly the Kronecker and plethysm coefficients, play a crucial role in Geometric Complexity Theory in the quest for separating algebraic complexity classes or simply separating two explicit polynomials. In order to achieve such separation we need to understand the multiplicities of irreducible components in the coordinate rings of the orbit closures of the given polynomials.
As it turned out just considering whether multiplicities are 0 or not is not enough in most cases of interest. This implies that we need to understand better what these multiplicities are and how large they can be.

One aspect of this understanding would be to find their combinatorial interpretation. For the Kronecker coefficients this has been an open problem in Algebraic Combinatorics and Representation Theory for more than 80 years. The fact that just deciding whether Kronecker coefficients is $\NP$-hard, and that the value of the character square of the symmetric group is not $\SP$, is evidence that sometimes these problems, as fundamental as they are, may not be doable the way we expect. Computational Complexity theory can help answer these questions and would be especially useful for \emph{negative} answers, if the situation happens to be such. 
 
Finally, moving beyond positivity and complexity of structure constants, in the lack of exact formulas, we turn towards their asymptotic properties and effective bounds. Estimating how large these multiplicities are for certain families is yet another big open problem, see~\cite{Pan22}. Such estimates could potentially close the cycle to GCT.


\begin{thebibliography}{Bou12300}

\bibitem[Aar16]{Aa}
S.~Aaronson, $\P\overset{?}=\NP$, in \emph{Open problems in mathematics},
Springer, Cham, 2016, 1--122.

\bibitem[AKS04]{AKS} M.~Agrawal, N.~Kayal and N.~Saxena,  PRIMES is in P. {\it Annals of mathematics} (2004), 781--793.



\bibitem[BO05]{BO} C.~Ballantine, R.~Orellana,  A combinatorial interpretation for the coefficients in the Kronecker product $s_{(n-p,p)}*s_\la$. \emph{S\'em. Lothar. Combin.} {\bf 54A} (2005/07), 29 pp.




\bibitem[B97]{B97} R.~Beals,  Quantum computation of Fourier transforms over symmetric groups, \emph{Proc. of the twenty-ninth annual ACM symposium on Theory of computing} (1997), pp. 48--53.

\bibitem[BB04]{BB} C.~Bessenrodt and C.~Behns,
On the Durfee size of Kronecker products of characters of
the symmetric group and its double covers,
{\em J.~Algebra}~{\bf 280} (2004), 132--144.

\bibitem[BB17]{BB17} 
C.~Bessenrodt, C.~Bowman, Multiplicity-free Kronecker products of characters of the symmetric groups. \emph{Adv. Math.} {\bf 322} (2017), pp. 473--529. 



\bibitem[Bla17]{Bla}
Jonah Blasiak.
\newblock {K}ronecker coefficients for one hook shape.
\newblock {\em S{\'e}m. Lothar. Combin}, 77:2016--2017, 2017.


\bibitem[BL18]{BL} J.~Blasiak, R.I.~Liu, Kronecker coefficients and noncommutative super Schur functions. \emph{ J. Combin. Theory Ser. A}{\bf 158} (2018), 315--361.

\bibitem[BMS15]{BMS}
J.~Blasiak, K.~Mulmuley, M.~Sohoni, Geometric complexity theory IV: nonstandard quantum group for the Kronecker problem. \emph{Mem. Amer. Math. Soc.}{\bf 235} (2015), no. 1109.


\bibitem[BI18]{BI}
M.~Bl\"aser and C.~Ikenmeyer,
\emph{Introduction to geometric complexity theory},
Summer school lecture notes, 2018, 148~pp.;
{\tt https://tinyurl.com/nhe2wxvw}



\bibitem[BCGHZ]{BCGHZ} S.~Bravyi, A.~Chowdhury,  D.~Gosset, V.~Havlicek, G.~Zhu, Quantum complexity of the Kronecker coefficients. (2023) \texttt{arXiv:2302.11454}.

\bibitem[B\"{u}r00a]{Bur00} P.~B\"urgisser,
\emph{Completeness and reduction in algebraic complexity theory},
Springer, Berlin, 2000, 168~pp.

\bibitem[B\"{u}r00b]{Bur0} P.~B\"urgisser, Cook's versus Valiant's
hypothesis, \emph{Theor.\ Comp.\ Sci.}~{\bf 235} (2000), 71--88.

\bibitem[B\"ur15]{Bur} P.~B\"urgisser,
Permanent versus determinant, obstructions, and Kronecker coefficients,
\emph{S\'em.\ Lothar.\ Combin.}~\textbf{75} (2015), Art.~B75a, 19 pp.

\bibitem[BCMW]{BCMW} P.~Burgisser, M.~Christandl, K.D.~Mulmuley, M.~Walter, Membership in moment polytopes is in NP and coNP. \emph{SIAM Journal on Computing}  (2017), {\bf 46}(3), 972-991.

\bibitem[BCS97]{BCS}
P.~B\"urgisser, M.~Clausen and M.~A.~Shokrollahi,
\emph{Algebraic complexity theory}, Springer, Berlin, 1997.

\bibitem[BI08]{BI08} P.~Bürgisser, C.~Ikenmeyer,  The complexity of computing Kronecker coefficients. \emph{ 20th Annual International Conference on Formal Power Series and Algebraic Combinatorics} (FPSAC 2008), pp. 357--368.



\bibitem[BIP19]{BIP}
P.~B\"urgisser, C.~Ikenmeyer and G.~Panova,
No occurrence obstructions in geometric complexity theory,
\emph{J.~AMS}~\textbf{32} (2019), 163--193; extended abstract in
\emph{57-th FOCS} (2016), IEEE, Los Alamitos, CA, 386--5.



\bibitem[BLMW11]{BLMW}
P.~B\"urgisser, J.~M.~Landsberg, L.~Manivel and J.~Weyman,
An overview of mathematical issues arising in the Geometric Complexity Theory
approach to $\VP\overset{?}=\VNP$, \emph{SIAM J.\ Comput.} {\bf 40} (2011), 1179--1209.


\bibitem[CHM07]{CHM}
 M.~Christandl, A.~Harrow, G.~Mitchison,  Nonzero Kronecker coefficients and what they tell us about spectra. \emph{Comm. Math. Phys.} {\bf 270 }(2007), no. 3, pp.575--585. 

\bibitem[CDW12]{CDW}
M.~Christandl, B.~Doran and M.~Walter, Computing Multiplicities of
Lie Group Representations, in \emph{Proc.\ 53-rd FOCS} (2012), IEEE, 639--648.

\bibitem[COSSZ]{COSSZ}
L.~Colmenarejo, R.~Orellana, F.~Saliola, .~Schilling, M.~Zabrocki,
The mystery of plethysm coefficients (2022), \texttt{arXiv:2208.07258} .



\bibitem[DM06]{DM06}
J.~A.~De~Loera and T.~B.~McAllister,
On the computation of Clebsch--Gordan coefficients and the dilation effect,
\emph{Experimental Math.}~\textbf{15} (2006), 7--19.






\bibitem[DIP19]{DIP}
J.~D\"orfler, C.~Ikenmeyer and G.~Panova,
{On geometric complexity theory: Multiplicity obstructions are stronger than occurrence obstructions} in
\emph{Proc.\ 46-th ICALP} (2019), Art.~51, 14~pp.


\bibitem[FI20]{FI20}
N.~Fischer and C.~Ikenmeyer,
The computational complexity of plethysm coefficients,
\emph{Comput.\ Complexity}~\textbf{29} (2020), no.~2,
Paper~8, 43~pp.

\bibitem[FRT54]{FRT}
J.~S.~Frame, G.~de~B. Robinson and R.~M. Thrall,
The hook graphs of the symmetric group,
{\em Canad.\ J.\ Math.}~\textbf{6} (1954), 316--324.



\bibitem[Ful97]{Ful}
W.~Fulton, \emph{Young tableaux},
Cambridge Univ.\ Press, Cambridge, UK, 1997, 260~pp.

\bibitem[Ful00]{Ful00}
W.~Fulton, Eigenvalues, invariant factors, highest weights, and Schubert calculus,
\emph{Bull.\ AMS}~\textbf{37} (2000), 209--249.




\bibitem[GJ79]{GJ}
M.~R.~Garey and D.~S.~Johnson, \emph{Computers and intractability},
Freeman, San Francisco, CA, 1979.

\bibitem[GR85]{GR85}
A.~Garsia and J.~Remmel,
\newblock Shuffles of permutations and the {K}ronecker product,
\newblock {\em Graphs and Combinatorics} (1985), 1:217--263.



\bibitem[GIP17]{GIP}
 F.~Gesmundo, C.~Ikenmeyer and G.~Panova,
{Geometric complexity theory and matrix powering},  \emph{Diff.\ Geom.\ Applications}~{\bf 55} (2017), 106--127.





\bibitem[Gre11]{Gre}
B.~Grenet, An upper bound for the permanent versus determinant problem,
manuscript (2011); \ts
{\tt http://www.lirmm.fr/~grenet/publis/Gre11.pdf}

\bibitem[HCW]{HCW} 
A.~Harrow, M.~Christandl, M.~Walter, Personal communication, \url{https://qi.ruhr-uni-bochum.de/talks/momenttalk_simons.pdf}  (2015).   

\bibitem[HSTZ13]{HSTZ}
G.~Heide, J.~Saxl, P.~H.~Tiep, A.~Zalesski,
Conjugacy action, induced representations and the Steinberg square for simple groups of Lie type, 
\emph{Proc. Lond. Math. Soc.} (3) {\bf 106} (2013), no. 4, 908--930.

\bibitem[I15]{I15}
C.~Ikenmeyer, The Saxl conjecture and the dominance order, \emph{Disc. Math.} (11) {\bf 6} (2015), pp.1970--1975.

\bibitem[IMW17]{IMW17}
C.~Ikenmeyer, K.~D.~Mulmuley and M.~Walter,
On vanishing of Kronecker coefficients,
\emph{Computational Complexity}~\textbf{26} (2017), 949--992.


\bibitem[IP22]{IP22}
C.~Ikenmeyer and I.~Pak,
What is in~$\SP$ and what is not?,  preprint (2022),
82~pp.; extended abstract to appear in \emph{Proc.\ 63rd FOCS} (2022); \ts {\tt arXiv:2204.13149}.

\bibitem[IPP22]{IPP22}
C.~Ikenmeyer, I.~Pak and G.~Panova,
Positivity of the symmetric group characters is as hard as the polynomial time hierarchy,
preprint (2022), 15 pp.; extended abstract to appear in \emph{Proc.\ 34th SODA} (2023); \ts {\tt arXiv:2207.05423}.

\bibitem[IP17]{IP}
C.~Ikenmeyer and G.~Panova, Rectangular Kronecker coefficients and plethysms
in geometric complexity theory, \emph{Adv.\ Math.}~\textbf{319} (2017), 40--66;
extended abstract in \emph{57th FOCS} (2016), IEEE, Los Alamitos, CA, 396--4055.

\bibitem[IS23]{IS}
C.~Ikenmeyer and S.~Subramanian, A remark on the quantum complexity of Kronecker coefficients, preprint (2023). 





\bibitem[KL14]{KL} H.~Kadish, J.~M.~Landsberg. Padded polynomials, their cousins, and geometric complexity theory. \emph{Comm. Algebra} {\bf 42}(5):2171--2180, (2014).

\bibitem[KT99]{KT99}
A.~Knutson and T.~Tao,
The honeycomb model of $\GL_n(\cc)$ tensor products~I:
Proof of the saturation conjecture,
{\em J.~AMS} {\bf 12} (1999), 1055--1090.



\bibitem[Las79]{Las79}
Alain Lascoux.
\newblock Produit de {K}ronecker des repr{\'e}sentations du groupe
  sym{\'e}trique.
\newblock In {\em S{\'e}minaire d'Alg{\`e}bre Paul Dubreil et Marie-Paule
  Malliavin: Proceedings, Paris 1979 (32{\`e}me Ann{\'e}e)}, pages 319--329.
  Springer, 1979.

\bibitem[Lit58]{Lit58}
Dudley~E. Littlewood.
\newblock Products and plethysms of characters with orthogonal, symplectic and
  symmetric groups.
\newblock {\em Canadian Journal of Mathematics}, 10:17--32, 1958.

\bibitem[LR34]{LR34} Littlewood, D. E., Richardson, A. R. (1934), Group Characters and Algebra, {\em Philosophical Transactions of the Royal Society of London  Series A}, {\bf 233} (721-730): 99-141.

\bibitem[Liu17]{Liu17}
Ricky~Ini Liu.
\newblock A simplified {K}ronecker rule for one hook shape.
\newblock {\em Proceedings of the American Mathematical Society},
  145(9):3657--3664, 2017.


\bibitem[LS17]{LS}
S.~Luo,M.~Sellke, The Saxl conjecture for fourth powers via the semigroup property. \emph{J. Algebraic Combin.} {\bf 45} (2017), no. 1, 33--80.



\bibitem[Mac95]{Mac}
I.~G.~Macdonald,
\emph{Symmetric functions and Hall polynomials} (Second edition),
Oxford University Press, New York, 1995.



\bibitem[MPP19a]{MPPe}
S.~Melczer, G.~Panova and R.~Pemantle,
{Counting partitions inside a rectangle},
\emph{SIAM J.\ Discrete Math.}~\textbf{34} (2020), 2388--2410;
extended abstact in {\em Proc.\ 31st FPSAC}, 2019.




\bibitem[MR04]{MR}
T.~Mignon and N.~Ressayre, A quadratic bound for the determinant and permanent problem,
\emph{Int.\ Math.\ Res.\ Notices}~\textbf{2004}, no.~79, 4241--4253.









\bibitem[Mul11]{Mul11}
K.~D.~Mulmuley, Geometric Complexity Theory VI: The flip via positivity,
preprint (2011), 40~pp., \ts available at \ts
{\tt http://gct.cs.uchicago.edu/gct6.pdf};
cf.\ {\tt arXiv:0704.0229}, 139~pp.


\bibitem[Mul17]{Mul} K.~Mulmuley,
Geometric Complexity Theory V.\ Efficient algorithms for Noether normalization,
\emph{J.\ AMS}~{\bf 30} (2017), 225--309.

\bibitem[MNS12]{MNS12}
K.~D.~Mulmuley, H.~Narayanan and M.~Sohoni,
Geometric complexity theory~III.  On deciding nonvanishing of a
Littlewood-Richardson coefficient, \emph{J.\ Algebraic Combin.}~\textbf{36}
(2012), 103--110.

\bibitem[MS01]{MS1}
K.~D.~Mulmuley and M.~Sohoni, Geometric complexity theory.~{I} \ts
{A}n approach to the {\P} vs.\ {\NP}  and related problems,
{\em SIAM J.\ Comput.}~\textbf{31}, 2001, 496--526.

\bibitem[MS08]{MS2}
K.~D.~Mulmuley and M.~Sohoni, Geometric complexity theory.~{II} \ts
{T}owards explicit obstructions for embeddings among class varieties,
{\em SIAM J.\ Comput.}~\textbf{38}, 2008, 1175--1206.


\bibitem[Mur38]{Mur38}
F.~D.~Murnaghan,
The analysis of the Kronecker product of irreducible representations of the symmetric group,
\emph{Amer.\ J.\ Math.}~\textbf{60} (1938), 761--784.

\bibitem[Mur56]{Mur56}
F.~D.~Murnaghan,
On the Kronecker product of irreducible representations of the symmetric group,
\emph{Proc.\ Natl.\ Acad.\ Sci.\ USA}~\textbf{42} (1956), 95--98.

\bibitem[Nar06]{Nar}
H.~Narayanan,
On the complexity of computing Kostka numbers and Littlewood--Richardson coefficients,
\emph{J.~Algebraic Combin.}~\textbf{24} (2006), 347--354.

\bibitem[Nis91]{Nis}
N.~Nisan,
Lower bounds for non-commutative computation, in
\emph{Proc.\ 23rd  STOC} (1991), ACM, 410--418.




\bibitem[Pak22+]{Pak22}
I.~Pak, What is a combinatorial interpretation?, preprint (2022), 58 pp.;
to appear in \emph{Proc.\ OPAC}, AMS, Providence, RI.

\bibitem[PP13]{PP13}
I.~Pak and G.~Panova, Strict unimodality of $q$-binomial coefficients,
\emph{C.\ts\/R.\/ Math.\/ Acad.\/ Sci.\/ Paris} \textbf{351} (2013), 415--418.

\bibitem[PP14]{PP14}
I.~Pak and G.~Panova,
Unimodality via Kronecker products, \emph{J.~Algebraic Combin.}~\textbf{40}
(2014), 1103--1120.

\bibitem[PP17a]{PPq} I.~Pak, G.~Panova, Bounds on certain classes of Kronecker and q-binomial coefficients,
\emph{Journal of Combinatorial Theory, Series A}{\bf 147} (2017), pp.1--17.

\bibitem[PP17b]{PP17}
I.~Pak and G.~Panova,
On the complexity of computing Kronecker coefficients,
\emph{Comput.\ Complexity}~\textbf{26} (2017), 1--36.

\bibitem[PP23]{PPd} I.~Pak, G.~Panova,  Durfee squares, symmetric partitions and bounds on Kronecker coefficients, \emph{J. of Algebra} (2023), {\bf 629} 358-380.

\bibitem[PPV16]{PPV} I.~Pak, G.~Panova, E.~Vallejo, Kronecker products, characters, partitions, and the tensor square conjectures,
\emph{Advances in Mathematics} {\bf 288} (2016), pp.702--731.



\bibitem[Pan15]{Pporto}
Greta Panova.
\newblock {K}ronecker coefficients: combinatorics, complexity and beyond.
Joint AMS-EMS meeting, Porto, Portugal (2015),
  \url{https://drive.google.com/file/d/1T2bVbLa4Bozy2_VBzT_Z0aezw8Y7ysrX/view}

\bibitem[Pan23]{Pan22} G.~Panova, Complexity and asymptotics of structure constants, in \emph{Open Problems in Algebraic Combinatorics Proceedings} (2023), \texttt{	arXiv:2305.02553}.  



\bibitem[Rem89]{R89}
J.~Remmel,
\newblock A formula for the {K}ronecker products of {S}chur functions of hook
  shapes,
\newblock {\em J. Algebra}(1989), {\bf 120}(1):100--118.


\bibitem[RW94]{RW} J.~Remmel, T.~Whitehead,  On the Kronecker product of Schur functions of two row shapes. \emph{Bull. Belg. Math. Soc. Simon Stevin} {\bf 1}(1994), no. 5, 649--683.

\bibitem[Sag01]{Sag}
B.~E.~Sagan, \emph{The symmetric group} (Second ed.), Springer, New York, 2001.


\bibitem[SW01]{SW}
A.~Shpilka and A.~Wigderson, Depth-$3$ arithmetic circuits over fields of
characteristic zero, \emph{Comput. Complexity}~\textbf{10} (2001), 1--27.

\bibitem[SY09]{SY}
A.~Shpilka and A.~Yehudayoff, Arithmetic Circuits: A survey of recent
results and open questions, \emph{Foundations and Trends in Theoretical
Computer Science}~\textbf{5} (2009), 207--388.



\bibitem[Sta99]{S1}
R.~P.~Stanley, {\em Enumerative Combinatorics}, vol.~1 (Second ed.)
and~2, Cambridge Univ.\ Press, 2012 and~1999.

\bibitem[Sta00]{Sta00}
R.~P.~Stanley,
Positivity problems and conjectures in algebraic combinatorics,
in \emph{Mathematics: frontiers and perspectivies},
AMS, Providence, RI, 2000, 295--319.

\bibitem[Tew15]{Tew15}
Vasu Tewari.
\newblock {K}ronecker coefficients for some near-rectangular partitions.
\newblock {\em Journal of Algebra}, 429:287--317, 2015.




\bibitem[Val79a]{V1}
L.~G.~Valiant, Completeness classes in algebra. \emph{Proc.\ 11th  STOC} (1979), 249--261.

\bibitem[Val79b]{V2}
L.~G.~Valiant,   The  complexity  of  computing  the  permanent.
\emph{Theor.\  Comp.\ Sci.}~{\bf 8} (1979), 189--201.

\bibitem[VSBR]{VSBR} L. G. Valiant, S. Skyum, S. Berkowitz, C. Rackoff. Fast Parallel Computation of
Polynomials Using Few Processors. \emph{SIAM Journal on Computing} {\bf 12}(4): 641-644 (1983).
%


\bibitem[Val97]{Val97}
E.~Vallejo, Reductions of additive sets, sets of uniqueness and pyramids,
\emph{Discrete Math.}~\textbf{173} (1997), 257--267.

\bibitem[Val99]{Val99}
E.~Vallejo,
Stability of Kronecker products of irreducible characters of
the symmetric group,
\emph{Electron.\ J.\ Combin.}~\textbf{6} (1999), RP~39, 7~pp.


\bibitem[VW14]{VW} Vergne, M.,  Walter, M. Inequalities for moment cones of finite-dimensional representations. (2014) \texttt{arXiv:1410.8144}.

\bibitem[Wig19]{Wig}
A.~Wigderson, \emph{Mathematics and Computation}, monograph draft, 2019.



\bibitem[Z89]{Z} D.~Zeilberger, A one-line high school algebra proof of the unimodality of the Gaussian polynomials $\binom{n}{k}_q$ for $k< 20$. (1989) In \emph{q-Series and Partitions } (pp. 67-72). Springer US.

\end{thebibliography}
\end{document}